\documentclass[a4paper,11pt]{article}
\usepackage{mathrsfs}
\usepackage{latexsym}
\usepackage{amsmath,amssymb}
\usepackage{amsthm}
\usepackage{amsfonts}
\usepackage[usenames]{color}
\usepackage{amssymb}
\usepackage{graphicx}
\usepackage{amsmath}
\usepackage{amsfonts}
\usepackage{amsthm}
\usepackage{mathrsfs}
\usepackage{dsfont}
\usepackage{indentfirst}
\usepackage[square, comma, sort&compress, numbers]{natbib}

\newcommand\D{\tilde{\Delta}}

\def\lan{\langle}
\def\ran{\rangle}
\def\p{\partial }

\def\phi{\varphi}

\def\D{\nabla}

\newtheoremstyle{mythm}{1.5ex plus 1ex minus .2ex}{1.5ex plus 1ex
minus .2ex}{\kai}{\parindent}{\song\bfseries}{}{1em}{}
\numberwithin{equation}{section}
\newtheorem{definition}{Definition}[section]
\newtheorem{theorem}{Theorem}[section]
\newtheorem{lemma}{Lemma}[section]

\newtheorem{proposition}{Proposition}[section]
\newtheorem{remark}{Remark} [section]
\newtheorem{corollary}{Corollary}[section]
\allowdisplaybreaks[4]

\begin{document}
\title
{Uniqueness of positive solutions for  m-Laplacian equations with polynomial non-linearity}
\author{Wei Ke}
\date{}
\maketitle

\begin{abstract}
 We consider the uniqueness of the following positive solutions of $m$-Laplacian equation:
\begin{equation}\label{eq1}
	\left\{
	\begin{aligned}
		-\Delta _m u&=\lambda u^{m-1}+u^{p-1} \quad \text{in} \quad  \Omega\\
		u&=0 \quad  \text{on} \quad  \partial \Omega  \\
	\end{aligned}
	\right.
\end{equation}
where $m>2$ is a constant. 
When $p\rightarrow m$, the uniqueness of positive solutions of $(\ref{eq1})$ is shown which is based on the essential uniqueness of first eigenfunction for $m$-Laplacian equation. Futhermore, we also prove the uniqueness results when $(\ref{eq1})$ is a perturbation of Laplacian equation.
\end{abstract}
\section{Introduction}
 In this paper, we study the uniqueness of positive solutions of $m$-Laplacian problem
\begin{equation}\label{eq2}
\left\{
\begin{aligned}
-\Delta _m u&=\lambda u^{m-1}+u^{p-1} \quad \text{in} \quad  \Omega\\
       	u&=0 \quad  \text{on} \quad  \partial \Omega  \\
\end{aligned}
\right.
\end{equation}
where $\Omega$ is a bounded $C^2$ domain in $\mathbb{R}^N$, $\Delta _m u:=\text{div}(|\nabla u|^{m-2}\nabla u)$.

A solution $u\in W^{1,m}_0(\Omega)$ of $(\ref{eq2})$ is called nondegenerate if the linearized problem
\begin{equation}\label{eq0}
	\left\{
	\begin{aligned}
		-\Delta _m v&=\lambda (m-1)u^{m-2}v+(p-1)u^{p-2}v \quad \text{in} \quad  \Omega\\
		v&=0 \quad  \text{on} \quad  \partial \Omega  \\
	\end{aligned}
	\right.
\end{equation}
admits only the trivial solution in $W^{1,m}_0(\Omega)$.

The uniqueness of non-trivial solutions of 
\begin{equation}
	\label{eq0011}
	\left\{
	\begin{aligned}
		-\Delta_mu&=f(u) \quad \text{in}  \quad \Omega,\\
		u&=0 \quad \text{on} \quad \partial \Omega
	\end{aligned}
	\right.
\end{equation}
remains widely open.
We denote by $\mathcal{D}^{1,m}_0(\Omega)$ the completion of $C^{\infty}_0(\Omega)$ with respect to the norm 
\[\phi\mapsto||\nabla \phi||_{L^m(\Omega)}.
\]
The space $\mathcal{D}^{1,m}_0(\Omega)$ is continuously embedded into $L^p(\Omega)$. 
There is a sharp constant for this embedding i.e., the quantity defined by 
\[\lambda_{m,p}(\Omega)=\inf_{u\in\mathcal{D}^{1,m}_0(\Omega)}\bigg\{\int_{\Omega}|\nabla u|^mdx:\int_{\Omega}|u|^pdx=1\bigg\}.
\]
Any minimizer of the sharp constant solves the following quasilinear version of the \textit{Lane-Emden} equation 
\begin{equation}
	\label{eq0012}
	\begin{aligned}
		-\Delta_mu=\lambda_{m,p}(\Omega)|u|^{p-2}u \quad \text{in}  \quad \Omega.
	\end{aligned}
\end{equation}

In Laplacian case, i.e. $m=2$, extremals are unique when $1\leq p<2$. One can infer uniqueness of positive solutions to the equation $(\ref{eq0012})$  by Brezis and Oswald in \cite{brezis1986remarks}. By this result and the fact that extremal for $\lambda_{2,p}$ never changes sign (see \citep [Proposition 2.3]{brasco2019overview}), we obtain the simplicity property. For $p=2$, $\lambda_{2,2}$ is nothing but the first eigenvalue of the Dirichlet-Laplacian. We can get the simplicity property by Linear Spectral Theory (see for example \cite [Theorem 1.2.5]{henrot2006extremum}). In super-homogeneous regime, uniqueness of extremal is known to hold in balls (see for example \citep[Theorem 2 \& Corollary 1]{kwong1992uniqueness}) and for planar convex sets (see for example \cite[Theorem 1]{lin1994uniqueness}). However, there are examples that extremal may not be unique. Nazarov constructs a counter-example in \cite[Proposition 1.2]{mikayelyan2015hopf}. Another counter-example can be found in a starshaped set consisting of two hypercubes overlapping in a small region near one corner in \cite[Example 4.7]{brasco2019overview}. 

For general $1<m<\infty$, uniqueness of minimizers holds for subhomogeneous case in $1\leq p<m$. This is a consequence of stronger uniqueness result for positive solutions of (\ref{eq0012}) contained in
 \cite[Theorem 4]{idogawa1995first} which is the quasilinear counterpart of the result by Brezis and Oswald. We can infer simplicity by combining this result with the fact that extremals for $\lambda_{m,p}$ (see for example the proof of \cite[Theorem 1]{idogawa1995first}). For $p=m$, the quantity $\lambda_{m,m}(\Omega)$ is the first eigenvalue of the $m$-Laplacian with Dirichlet conditions, then we can get the simplicity. In super-homogeneous regime, simplicity holds in a ball, but fails in general. The counter-example by Nazarov works in this case, too. We also refer to Kawohl's paper \cite{kawohl2000symmetry} for the very same example. 
 
For fractional Laplacian equation, Dieb, Ianni and Salda\~na discussed positive solutions of the Dirichlet Lane-Emden-type fractional problem
\begin{equation}\label{eq004}
	\left\{
	\begin{aligned}
		(-\Delta )^s  u+\lambda u&=u^p \quad \text{in} \quad  \Omega,\\
		u&=0 \quad  \text{on} \quad  \mathbb{R}^N \backslash \Omega,  \\
	\end{aligned}
	\right.
\end{equation}
where $s\in (0,1)$, $p>1$, and $\lambda\in \mathbb{R}$, in \cite{dieb2022uniqueness}. They obatin the uniqueness and nondegeneracy hold for the asymptotically linear problem. Futhermore, they prove that the uniqueness and nondegeneracy results when the fractional parameter sufficiently close to 1. But the uniqueness of solutions of $(\ref{eq004})$ does not hold in general, see \cite{angeles2023small} and \cite{long2019critical}.

For fully nonlinear case, 
\begin{equation}
	\label{eq0013}
	\left\{
	\begin{aligned}
		\det D^2u&=f(-u) \quad \text{in}  \quad \Omega,\\
		u&=0 \quad \text{on} \quad \partial \Omega,
	\end{aligned}
	\right.
\end{equation}
there are few results. The uniqueness of solutions of $(\ref{eq0013})$ is known for $f=\Lambda_qt^q $  and $0\le q<N$ in  \cite{Tso1990}. When $q=N$, this is the eigenvalue problem of Monge-Ampère equation which was studied in \cite{lions1985two}.
For super-linear case, Huang studied the uniqueness  results on smooth uniformly convex domain in \cite{huang2019uniqueness}. Cheng, Huang, and Xu obtained the uniqueness results on symmetric convex domain in $\mathbb{R}^2$ in \cite{2023Huang}.

\noindent 1.1. \textbf{Main results.}
Let $\lambda_{m}(\Omega)>0$ denote the first Dirichlet eigenvalue for $-\Delta_m$ in $\Omega$. For $\lambda \in (0,\lambda_m)$ and $p\in (1,\frac{Nm}{N-m})$, the equation $(\ref{eq2})$ has a non-trivial solution in \cite[Section 6]{garcia1987existence}. In this paper, we study $m$-Laplacian equation and obtain two following theorems. 
\begin{theorem}\label{thm1}
	Let $N\geq 2$, $p\in (2,\frac{2N}{N-2})$, $\lambda\in (0,\lambda_2(\Omega))$, and $\Omega$ be such that the problem
	\begin{equation}\label{eq04}
		\left\{
		\begin{aligned}
			-\Delta  u&=\lambda u+u^{p-1} \quad \text{in} \quad  \Omega\\
			u&=0 \quad  \text{on} \quad  \partial \Omega  \\
		\end{aligned}
		\right.
	\end{equation}
	has a unique positive solution which is nondegenerate. Then, there is $\delta_1 =\delta_1 (\Omega,\lambda,p)\in (0,1)$ such that, for $m\in (2,2+\delta_1]$, the problem $(\ref{eq2})$ has a unique positive solution.
\end{theorem}
\begin{theorem}\label{thm2}
	Let $m\in (2,+\infty)$, $N\geq 2$, and $\Omega$ be a bounded domain of class $C^2$, and $\lambda \in (0,\lambda_m(\Omega))$. There is a $p_0=p_0(\Omega,m,\lambda)>m$ such that problem $(\ref{eq2})$ has a unique positive solution for every $p\in (m,p_0)$.
\end{theorem}

\noindent 1.2. \textbf{Some comments on the proof.} The proof of Theorem 1.2 is largerly inspired by \cite[Lemma 3]{lin1994uniqueness}, dealing with the case of Laplacian. Brasco and Lindgren use the method to deal with $m$-Laplacian in \cite{brasco2023uniqueness}. 

By assuming that simplicity fails, there exists a sequence $\{q_n\}_{n\in \mathbb{N}}$ such that $q_n\rightarrow q$ and $\lambda_{p,q}$ admits two linearly independent positive minimizers $u_n$ and $v_n$. Let $\tilde{u}_n=\frac{u_n}{M_n}$ $\tilde{v}_n=\frac{v_n}{M_n}$, where $M_n$ is $L^{\infty}$ norm of $u_n$ . Their difference $\tilde{u}_n-\tilde{v}_n$ must be sign-changing on the domain and solves the linear equation 
\begin{equation}
	\begin{split}
		-\text{div} (A_n\nabla (\tilde{u}_n-\tilde{v}_n))&=(m-1)\lambda \int_{0}^{1}z_n^{m-2}(\tilde{u}_n-\tilde{v}_n)dt\\
		&+M_n^{p_n-m}(p_n-1)\int_{0}^{1}z_n^{p_n-2}(\tilde{u}_n-\tilde{v}_n)dt,
		\nonumber
	\end{split}
\end{equation}
where 
\[z_n:=(1-t)\tilde{v}_n+t\tilde{u}_n.
\]
The coefficient matrix 
\[A_n=\int_{0}^{1}|\nabla z_n|^{m-2}\text{Id}+(m-2)|\nabla z_n|^{m-4} \nabla z_n \otimes \nabla z_ndt
\]
 is degenerate elliptic. More precisely 
\[\frac{1}{C}(|\nabla \tilde{u}_n|^{p-2}+|\nabla \tilde{v}_n|^{p-2})|\xi|^2\le \left \langle A_n \xi,\xi \right \rangle \le C(|\nabla \tilde{u}_n|^{p-2}+|\nabla \tilde{v}_n|^{p-2})|\xi|^2, \quad \text{for every } \xi \in \mathbb{R}^N.
\]
 This permits to obtain a compact embedding in some $L^t$ space, for weighted Sobolev spaces of functions such that
\[\int_{\Omega} (|\nabla \tilde{u}_n|^{p-2}+|\nabla \tilde{v}_n|^{p-2})|\nabla \phi|^2dx<+\infty.
\]
 It is possible to pass the limit in the linearized equation. The conclusion of the proof is similar to that of \cite{lin1994uniqueness}: obtain convergence to a  non-trivial limit function $\varphi$, which by construction is a sign-changing first eigenvalue of a weighted linear eigenvalue problem. 

In our paper, it is diffcult for us to build a uniform boundary estimate as the domain grows. We can build a uniform estimate in local region in next section. Futhermore, we do not have a uniform priori bound for all positive solutions. We need to learn from the method in \cite[Section 4]{dieb2022uniqueness} and \cite {lin1994uniqueness}. Let maximum value of the solution divides this equation. Then we can obtain the regularity in \cite{fan2007global}. 

\noindent1.3. \textbf{Plan of the paper.} In Section 2, we prove some preliminary facts. In particular, we devote this section to proving weighted embeddings. Here we need to go through the proofs of \cite{brasco2023uniqueness}  and \cite{damascelli2004}. In Section 3, we analyze the first eigenvalue of a weighted linear eigenvalue problem. This is a crucial ingredient for the proof of our main result. We prove the existence of the first eigenfunction and its uniqueness. In Section 4, we prove Theorem 1.2. In Section 5, we prove Theorem 1.1. Finally, we include appendice, contains the crucial regularity results  in \cite{brasco2023uniqueness}  and \cite{damascelli2004}.

\noindent\textbf{Acknowledgements.} The author is grateful to professor Genggeng Huang for suggesting this problem
and for helpful discussions.

\section{ Preliminaries}
In this section, we introduce some notations and preliminary results needed throughout this paper.

Let $\Omega$ be a bounded domain of class $C^2$ in $\mathbb{R} ^N$, $N\geq 2$, $2<m<\infty$, $p>2$. Denote  $x=(x_1,\cdots, x_{N-1},x_N)=(x',x_N)$ for $x\in \mathbb{R}^N$. 

Set 
\[u_+:=\max\{u,0\} \quad \text{and} \quad  u_-:=\max\{-u,0\}.
\]

We have the following local gradient estimate for m-laplacian equation. For the case of Laplacian, i.e., $m=2$, we can see Qing Han's book in \cite[Lemma 1.1.11 \& Theorem 1.1.14]{han2016nonlinear}.
\begin{lemma}\label{lem1}
   Let $\Omega\subset \mathbb{R}^N$ be a domain satisfying a uniform exterior sphere condition, namely, there is $r_0>0$ such that, for every $x^0\in \partial \Omega$, there is $y^0\in \mathbb{R}^N \backslash \Omega$ with $\overline{B_{r_0}(y^0)}\cap \overline{\Omega}={x^0}$. Let $u\in  C(\overline{\Omega})\cap C^{1,\alpha}(\Omega)$ be a solution of
   \begin{equation}\label{eq5}
   	\left\{
   	\begin{aligned}
   		-\Delta _m u&=f \quad \text{in} \quad  \Omega\\
   		u&=0 \quad  \text{on} \quad  \partial \Omega  \\
   	\end{aligned}
   	\right.
   \end{equation}
with $f\in L^{\infty}(\Omega)\cap C(\overline{\Omega})$, $\alpha\in(0,1)$. Let $M>0$ satisfy that $||u||_{L^{\infty}(\Omega)}+||f||_{L^{\infty}(\Omega)}<M$. Then for any $x\in \Omega\cap B_1(x^0)$, $|u(x)-u(x^0)|\leq C|x-x^0|$, where $C$ is a positive constant depending only on $m$, $N$, $r_0$ and $M$.
\end{lemma}
\noindent \emph{Proof.}  For the given $x^0\in \partial \Omega$, let $d(x)$ be the distance from $x$  to $\partial  B_{r_0}(y^0)$; i.e., $d(x):=|x-y^{0}|-r_0$. 

Claim: there exists a  $C^2$ function  $\omega =\psi(d)$ defined in $[0,1)$, with $\psi (0)=0$ , $\psi' >0$ in $(0,1)$ and $\Delta _m\psi <-C(m,N,r_0)<0$.

Proof of the Claim: 
\[\partial_i d(x)=\frac{x_i-y^0_{i}}{|x-y^0|},
\] 
\[\partial_{ij} d(x)=\frac{\delta_{ij}}{|x-y^0|}-\frac{(x_i-y^0_{i})(x_j-y^0_{j})}{|x-y^0|^3}.
\]
Then
\[\Delta_m \omega =\text{div}(|\psi '|^{m-1}\nabla d)=(m-1)(\psi')^{m-2}\left(\psi ''+\frac{N-1}{m-1}\frac{\psi '}{r_0+d}\right).
\]
Set $a=\frac{N-1}{(m-1)r_0}$, then we get 
\[\Delta _m\omega \le (m-1)(\psi')^{m-2}(\psi ''+a\psi'+b)-(m-1)(\psi')^{m-2}b
\]
where $b$ is a  positive constant.
We  choose $b$ and find a function $\psi$ in  $[0,1)$ such that
\[\psi''+a\psi'+b=0\quad \text{in} \quad (0,1).
\]
First, a particular solution of the ordinary differential equation above is given by
\[\psi(d)=-\frac{b}{a}d+\frac{A}{a}\left(1-e^{-ad}\right)
\]
for some constant $A$. Let $A=\frac{2b}{a}e^{a}$, then
$\psi$ satisfies all the requirements we imposed.
We get 
\begin{equation}
	\left\{
	\begin{aligned}
		\label{eq5.0}
		\Delta_m(C\psi)&\leq\Delta_m  u \quad \text{in} \quad \Omega\cap B_1(x^0),\\
		C\psi&\geq  u \quad \text{on} \quad \partial (\Omega \cap B_1(x^0)),
	\end{aligned}
\right.
\end{equation}
where $C$ is a positive constant large enough which depends only on $m$, $N$, $r_0$, and $M$.  Choose the test function $(u-C\psi)_+$ on the first equation of $(\ref{eq5.0})$. Let
\[A_n(x)=\int_{0}^{1}\nabla ^2\mathcal{H}(t\nabla (C\psi)+(1-t)\nabla u)dt,
\]
where $\mathcal{H}$ is defined in section 3 for $p=m$.
Thus 
\[\int_{\Omega \cap B_1(x^0)}A_n(x)(\nabla (C\psi)-\nabla u)\cdot \nabla (u-C\psi)_+\ge 0.
\]
Then  
\[\int_{\Omega \cap B_1(x^0)}\frac{1}{4^{m-1}}|\nabla (C\psi)|^{m-2}|\nabla ((u-C\psi)_+)|^2\le 0.
\]
Thus 
\[\nabla (u-C\psi)_+=0 \quad \text{a.e. in} \quad \Omega \cap B_1(x^0).
\]
We get 
\[u\leq C\psi  \quad \text{in} \quad \Omega \cap B_1(x^0).
\]
Similarly 
\[-u\leq C\psi  \quad \text{in} \quad \Omega \cap B_1(x^0).
\]
We obtain 
\[|u-u(x^0)|\leq C|\psi-\psi (x^0)|.
\]
This implies the desired result. \hfill $ \Box$
\begin{remark}
	\label{rem2}
	Let $\Omega$ be a domain satisfying uniform exterior sphere condition of $r_0$. Note that, if $\mu>1$, then the domain $\mu\Omega$ satisfies in particular the uniform exterior sphere condition with the same $r_0$.
\end{remark}
\begin{remark}
	\label{rem2.5}
	$K$ is a compact subset in $(1,+\infty)$. $\forall m \in K$, the constant $C$ in Lemma \ref{lem1} remains bounded uniformly. 
\end{remark}
In  \cite{fan2007global}, Fan considered the equation  
\begin{equation}
	\label{eq4}
	\begin{aligned}
			-\text{div} A(x,u,Du)=B(x,u,Du)\quad\text{in} \quad \Omega\\
	\end{aligned}
\end{equation}
with Dirichlet boundary condition
\begin{equation}
	\label{eq7}
	u=g   \quad \text{on} \quad \partial \Omega.
\end{equation}
 $p(x):\overline{\Omega}\rightarrow \mathbb{R}$ satisfies the condition 
\begin{equation}
	\label{eq8}
	1<p_{-}=p_-(\Omega):=\inf _{x\in \Omega} p(x)\leq p_{+}(\Omega):=\sup_{x\in \Omega}p(x)<\infty.
\end{equation}
The symbols $\alpha$, $\beta$, $\gamma$ or $\alpha_i$, $\beta_i$, $\gamma_i$ denote positive constants in $(0,1)$.

\noindent \textbf{Assumption ($P_H$).} $p$ is H\"older continuous on $\overline{\Omega}$, which is denoted by $p\in C^{0,\beta_1}(\overline{\Omega})$, that is, there exist a positive constant $L$ and an exponent $\beta_1\in (0,1)$ such that 
\[|p(x_1)-p(x_2)|\leq L|x_1-x_2|^{\beta_1} \quad \text{for} \quad x_1,x_2\in \overline{\Omega}
\]
\noindent \textbf{Assumption ($A^k$).} $A=(A_1,A_2,\cdots,A_N)\in C(\overline{\Omega} \times\mathbb{R}\times\mathbb{R}^N, \mathbb{R}^N)$. For every $(x,u)\in \overline{\Omega}\times \mathbb{R}$, $A(x,u,\cdot)\in C^1(\mathbb{R}^N\backslash\{0\}, \mathbb{R}^N)$, and there exist a nonnegative constant $k\geq 0$, a non-increasing continuous function $\lambda:[0,+\infty)\rightarrow (0,+\infty)$ and a non-decreasing continuous function $\Lambda:[0,+\infty)\rightarrow (0,+\infty)$ such that for all $x,x_1,x_2\in \overline{\Omega},u,u_1,u_2\in \mathbb{R}^N, \eta \in \mathbb{R}^N\backslash\{0\}$ and $\xi=(\xi_1,\xi_2,\cdots,\xi_N)\in \mathbb{R}^N,$ the following conditions are satisfied:
\begin{equation}
	A(x,u,0)=0
\end{equation}
\begin{equation}
	\sum_{i,j}\frac{\partial A_j}{\partial \eta_i}(x,u,\eta)\xi_i\xi_j\geq\lambda(|u|)(k+\eta ^2)^{\frac{p(x)-2}{2}}|\xi|^2,
\end{equation}
\begin{equation}
	\sum_{i,j}|\frac{\partial A_j}{\partial \eta_i}(x,u,\eta)|\leq\Lambda(|u|)(k+\eta ^2)^{\frac{p(x)-2}{2}},
\end{equation}
\begin{equation}
	\begin{split}
		|A&(x_1,u_1,\eta)-A(x_2,u_2,\eta)|\\
		&\leq \Lambda(\text{max}\{|u_1|,|u_2|\})(|x_1-x_2|^{\beta_1}+|u_1-u_2|^{\beta_2})\\
		&\times[(k+|\eta|^2)^{\frac{p(x_1)-2}{2}}+(k+|\eta|^2)^{\frac{p(x_2)-2}{2}}]|\eta|(1+|\text{log}(k+|\eta|^2)|).
	\end{split}
\end{equation}
\noindent \textbf{Assumption ($B$).} $B:\overline{\Omega}\times \mathbb{R} \times \mathbb{R}^N\rightarrow \mathbb{R},$ the function $B(x,u,\eta)$ is measureable in $x$ and is continuous in $(u,\eta)$, and 
\begin{equation}
	\label{eq9}
	|B(x,u,\eta)|\leq \Lambda (|u|)(1+|\eta|^{p(x)}), \quad \forall(x,u,\eta)\in \overline{\Omega}\times \mathbb{R}\times \mathbb{R}^N,
\end{equation}
where $\Lambda$ is as in the assumption $(A^k)$.

\begin{theorem}\label{thm3}
	\cite[Theorem 1.1]{fan2007global} Under Assumptions $(p_H), (A^k),$ and $(B)$, suppose that $u\in W^{1,p(x)}(\Omega)\cap L^{\infty}(\Omega)$ is a bounded generalized solution of $(\ref{eq4})$ and 
	\begin{equation}
		\label{eq10}
		\sup_{\Omega}|u(x)|:=\text{ess} \sup_{\Omega}|u(x)|\leq M \quad \text{for some} \quad M>0.
	\end{equation}
 Then $u\in C_{\mathrm{loc}}^{1,\alpha}(\Omega)$, where the H\"older  exponent $\alpha$ depends only on $p_-$, 
 $p_+$, $N$, $\lambda(M)$, $M$, $\Lambda(M)$, $L$, $\beta_1$ and $\beta_2$, and for given $\Omega_0 \Subset \Omega$, $|u|_{C^{1,\alpha}(\overline{\Omega}_0)}$- norm  of $u$, depends on $p_-$, $p_+$, $N$, $\lambda(M)$, $M$, $\Lambda(M)$, $L$, $\beta_1$, $\beta_2$ and $\text{dist}(\Omega_0,\partial \Omega)$.
\end{theorem}
\begin{theorem}\label{thm4}
	\cite[Theorem 1.2]{fan2007global} Under Assumptions $(p_H), (A^k),$ and $(B)$, also let the boundary $\partial \Omega$ of $\Omega$ be of class $C^{1,\gamma}$, and $g\in C^{1,\gamma } (\partial \Omega)$. Suppose that $u\in W^{1,p(x)}(\Omega)\cap L^{\infty}(\Omega)$ is a bounded generalized solution of $(\ref{eq4})$ and satisfies $(\ref{eq10})$. Then $u\in C^{1,\alpha}(\overline{\Omega})$, where the H\"older  exponent $\alpha$ and $|u|_{C^{1,\alpha}(\overline{\Omega})}$ depends only on $p_-$, $p_+$, $N$, $\lambda(M)$, $M$, $\Lambda(M)$, $L$, $\beta_1$, $\beta_2$,$\gamma$, $|g|_{C^{1,\gamma}(\partial \Omega)}$ and $\Omega$.
\end{theorem}
\begin{remark}
	\label{re4}
	m-Laplacian satisfies the Assumption $(p_H)$ and $(A^k)$.
\end{remark}
The $L^{\infty}$ norm of the solution for $(\ref{eq2})$ in $W^{1,m}_0(\Omega)$ can be found in \cite[Corollary 1.1]{guedda1989quasilinear}. Then, we can obtain the global $C^{1,\alpha}$ regularity by Theorem $\ref{thm4}$.

The following two lemmas are inspired in \cite{lin1994uniqueness}, dealing with the case of Laplacian.
\begin{lemma}
	\label{lem2} 
	Let $m>2$, $u_p\in W^{1,m}_0(\Omega)$ be a positive solution of $(\ref{eq2})$. Let $M_p=||u_p||_{L^{\infty}(\Omega)}$. Then, there is a $\delta_0=\delta_0(m,N,\Omega)>0$, such that $M_p^{p-m}\leq C$ for $ p\in [m,m+\delta_0]$, where $C= C(m, N,\Omega,\delta_0)>0$.
\end{lemma}
\noindent \emph{Proof.} By contradiction, assume that $ \{p_k\}\rightarrow m$, $M_{p_k}^{p_k-m}=||u_{p_k}||_{L^{\infty}(\Omega)}^{p_k-m}\rightarrow +\infty$, where $u_k$ is a sequence of positive solutions of $(\ref{eq2})$ with $p=p_k$. Let $\{x^k\}_{k\in \mathbb{N}}$ be a sequence in $\Omega$ such that $u_k(x^k)=M_{p_k}$. Let 
\[\tilde{u}_k(y)=\frac{u_k\left(\mu_ky+x^k\right)}{M_{p_k}}.
\]
Then $\tilde{u}_k$ is a function satisfying $0\leq\tilde{u}_k\leq1$, $\tilde{u}_k(0)=1$ and 
\begin{equation}
	\label{eq13}
    \left\{
    \begin{split}
    	-\Delta_m \tilde{u}_k&=\lambda \mu^m_k\tilde{u}^{m-1}_k+M_{p_k}^{p_k-m}\mu^m_k\tilde{u}^{p_k-1}_k:=f_k  \quad \text{in} \quad \Omega_k,\\
    	\tilde{u}_k&=0 \quad \text{on} \quad \partial \Omega_k,
    \end{split}
    \right.
\end{equation}
where $\Omega_k=\{y\in \mathbb{R}^N:\mu_ky+x^k\in \Omega\}. $
Let $M_{p_k}^{p_k-m}\mu^m_k=1$, then $\mu_k\rightarrow0 (\text{as} \quad k\rightarrow+\infty)$. $\tilde{u}_k$ satisfies
\begin{equation}
	\label{eq14}
	\left\{
	\begin{split}
		-\Delta_m \tilde{u}_k&=\lambda \mu^m_k\tilde{u}^{m-1}_k+\tilde{u}^{p_k-1}_k  \quad \text{in} \quad \Omega_k,\\
		\tilde{u}_k&=0 \quad \text{on} \quad \partial \Omega_k.
	\end{split}
	\right.
\end{equation}
Up to some subsequences, two situations may occur:
\[\text{either} \quad \text{dist}(x^k,\partial \Omega)\mu^{-1}_k\rightarrow +\infty \quad \text{or} \quad \text{dist}(x^k,\partial \Omega)\mu^{-1}_k\rightarrow d \geq 0.
\]

Case 1: $\lim_{k\rightarrow +\infty}\text{dist}(x^k,\partial \Omega)\mu^{-1}_k\rightarrow +\infty$. Then $\Omega_k\rightarrow \mathbb{R}^N$ as $k\rightarrow +\infty$
for any $r>1$ with $B_r(0)\subset\Omega_k$, $f_k$ satisfies $(\ref{eq9})$. By applying Theorem $\ref{thm3}$, there exist an $\alpha\in (0,1)$ and a constant $C=C(r)>0$ which don't depend on $k$  such that 
\[||\tilde{u}_k||_{C^{1,\alpha}(B_r(0))}\le C.
\]
By choosing a subsequence, $\tilde{u}_k\rightarrow \tilde{u}$ in $C^{1,\beta}_{\mathrm{loc}}(\mathbb{R}^N)$ for some $\beta\in (0,\alpha)$.
Moreover, $\tilde{u}(0)=1$, $\tilde{u}\ge 0$ in $\mathbb{R}^N$, $f_k \rightarrow \tilde{u}^{m-1}$ in $C^{1,\beta}_{\mathrm{loc}}(\mathbb{R}^N)$, and 
\[\int_{\mathbb{R}^N}|\nabla \tilde{u}|^{m-2}\nabla \tilde{u} \nabla \psi dx=\int_{\mathbb{R}^N}\tilde{u}^{m-1}\psi dx \quad \text{for all} \quad \psi \in C^1_c(\mathbb{R}^N),
\]
namely, $\tilde{u}$ is a distributional solution of
\begin{equation}
	\label{eq17}
	-\Delta_m \tilde{u}=\tilde{u}^{m-1} \quad \text{in} \quad \mathbb{R}^N.
\end{equation} 
We claimed that $\tilde{u}=0$. Suppose to the contrary that $\tilde{u}\not\equiv0$. Let $\varphi _R$ be the first eigenfunction of $-\Delta_m$ on $B_R(0)$, namely , $\varphi_R$ satisfies 
\begin{equation}
	\label{eq16}
	\left\{
	\begin{split}
		-\Delta_m \varphi_R&=\lambda_R \varphi^{m-1}_R \quad \text{in} \quad B_R(0),\\
		\varphi_R&=0 \quad \text{on} \quad \partial B_R(0),
	\end{split}
	\right.
\end{equation}
where $\lambda_R$ is the first eigenvalue of $-\Delta_m$ on $B_R(0)$ and $\lambda_R\rightarrow 0$ as $R\rightarrow +\infty$.
Testing $\psi_R=\varphi_R^m\tilde{u}^{1-m}$ on $(\ref{eq17})$, we obtain
\[	\int_{B_R(0)}|\nabla \tilde{u}|^{m-2}\nabla \tilde{u}\cdot \nabla(\varphi_R^m\tilde{u}^{1-m})dx= \int_{B_R(0)}\varphi^m_R dx.
\]
The left hand side denotes by LHS. Then, by Young's inequality,
\begin{equation}
	\label{eq21}
	\begin{split}
		\text{LHS}&=(1-m)\int_{B_R(0)}\varphi_R^m|\nabla \tilde{u}|^m \tilde{u} ^{-m}dx+m\int_{B_R(0)}|\nabla \tilde u|^{m-2}\tilde{u}^{1-m}\varphi^{m-1}_R\nabla \tilde{u}\cdot\nabla \varphi_Rdx,\\
		&\le (1-m)\int_{B_R(0)}\varphi_R^m|\nabla \tilde{u}|^m \tilde{u} ^{-m}dx +(m-1)\int_{B_R(0)}\varphi_R^m|\nabla \tilde{u}|^m \tilde{u} ^{-m}dx+C_m\int_{B_R(0)}|\nabla \varphi_R|^mdx.
	\end{split}
\end{equation}
We obtain
\[\int_{B_R(0)}|\nabla \varphi_R|^mdx\ge C_m^{-1}\int_{B_R(0)}|\varphi_R|^mdx.
\]
By $(\ref{eq16})$, 
\[\lambda_R\int_{B_R(0)}|\varphi_R|^mdx \ge C_m^{-1}\int_{B_R(0)}|\varphi_R|^mdx
\]
We let $R\rightarrow \infty$, and get the contradiction. Thus $\tilde{u}\equiv 0$.
But this contradicts  $\tilde{u}(0)=1$.

Case 2: $\text{dist}(x^k,\partial \Omega)\mu_k^{-1}\rightarrow d\ge 0$. Then $x^k\rightarrow x^0\in \partial \Omega$. Let $d_k$ be the distance from $x^k$ to $\partial \Omega$ and $d_k=\text{dist}(x^k,\tilde{x}^k)$, where $\tilde{x}^k\in \partial \Omega$. After translation, we may assume $\tilde{x}^k=0$, $\nu(0)=-e_N$, where $\nu(0)$ is the outward unit normal vector at $0$. 

Claim: $d>0$.
By contradiction, assume that $d=0$. We define the function $v_k(y)=u_k(\mu_ky)/M_{p_k}$ for $y\in \tilde{\Omega}_k$, where  $\tilde{\Omega}_k=\{y\in \mathbb{R}^N|\mu_ky\in \Omega\}$.Then $v_k$ satisfies that 
\begin{equation}
	\label{eq23}
	\left\{
	\begin{split}
		-\Delta_m v_k&=\lambda \mu^m_kv^{m-1}_k+v^{p_k-1}_k  \quad \text{in} \quad \tilde{\Omega}_k,\\
		v_k&=0 \quad \text{on} \quad \partial \tilde{\Omega}_k.
	\end{split}
	\right.
\end{equation}
Setting $z^k=\frac{x^k}{\mu_k}$, we obtain $v_k(z_k)=1$ and 
\[z^k\rightarrow 0 \quad \text{as} \quad k\rightarrow +\infty.
\]
 Near the origin, $\partial \Omega$ can be represented by $x_N=\rho(x')$, $\rho (0)=0$ and $\nabla \rho(0)=0$. 
We have $\mu_ky_N=\rho(\mu_ky')$. We obtain $y_N=\frac{\rho(\mu_ky')}{\mu_k} \rightarrow 0 $ locally uniformly as  $k\rightarrow +\infty$. By Lemma $\ref{lem1}$, Remark $\ref{rem2}$, for any $y\in \tilde{\Omega}_k\cap B_1(0)$,
$|v_k(y)-v_k(0)|\le C|y|$, where $C$ doesn't depend on $k$. Then,
\[|v_k(z^k)-v_k(0)|\le C|z^k|\rightarrow  0 \quad \text{as} \quad k \rightarrow +\infty.
\]
But this contradicts $v_k(z^k)=1$.

By choosing a subsequence, $v_k\rightarrow v$ locally uniformly $C^{1,\beta}(\{x_N>-d\})$, where $\mathbb{R}^N_{d+}$ denotes $\{x\in \mathbb{R}^N|x_N>-d\}$. By a translation, $v$ satisfies  
\begin{equation}
	\label{eq17.5}
	\left\{
	\begin{split}
		-\Delta_m v&=v^{m-1} \quad \text{in} \quad \mathbb{R}^N_+,\\
		v&=0 \quad \text{on} \quad x^N=0,\\
		v(de_N)&=1,\\
		v&> 0 \quad \text{in} \quad  \mathbb{R}^N_+.
	\end{split}
	\right.
\end{equation}
Let $\varphi_R$ be the first positive eigenfunction of $-\Delta_m$ on $B_R(Re_N)$. Then a similar argument as in Case 1 yields that $v\equiv 0$.
But this contradicts  $v(de_N)=1$.
The lemma is proved. \hfill $\Box$
\begin{lemma}
	\label{lem3}
	Let $m>2, \{p_k\}_{k\in \mathbb{N}}\subset [m,m+\delta_0]$ be a sequence such that $p_k\rightarrow m$, where  $\delta_0$ is defined in Lemma $\ref{lem2}$. Let $u_k$ be a sequence positive solutions of  $(\ref{eq2})$ with $p=p_k$. Let $M_{p_k}:=||u_k||_{L^{\infty}(\Omega)}$. Then exists a subsequence denote by $\{u_k\}_{k\in \mathbb{N}}$ such that $\tilde{u}_k:=\frac{u_k}{M_{p_k}}\rightarrow \tilde{u}$ in  $C^{1,\beta}(\Omega)$ for some $\beta\in (0,1)$, where $\tilde{u}$ is the first eigenfunction of $-\Delta_m$.
\end{lemma}
\noindent \emph{Proof.} $\tilde{u}_k$ satisfies 
\begin{equation}
	\label{eq 25}
	\left\{
	\begin{split}
		-\Delta_m \tilde{u}_k&=\lambda \tilde{u}_k^{m-1}+M^{p_k-m}_{p_k}\tilde{u}^{p_k-1}_k:=f_k \quad \text{in} \quad  \Omega,\\
		\tilde{u}_k&=0 \quad \text{on} \quad \partial \Omega.
	\end{split} 
	\right.
\end{equation}
$f_k$ has a uniform bound and satisfies $(\ref{eq9})$. By Theorem $\ref{thm4}$, $||\tilde{u}_k||_{C^{1,\alpha}(\Omega)}\le C$, where  $\alpha\in (0,1)$ and $C$ depend on $m$, $\delta_0$, $N$. Thus there exist $\beta\in (0,\alpha)$ and a subsequence $\tilde{u}_k\rightarrow \tilde{u}\ge 0$ in $C^{1,\beta}(\overline{\Omega})$ and $||\tilde{u}||_{L^{\infty}(\Omega)}=1$. By Lemma \ref{lem2}, $\tilde{u}$ satisfies 
\begin{equation}
	\label{eq 26}
	\left\{
	\begin{split}
		-\Delta_m \tilde{u}&=\lambda \tilde{u}^{m-1}+C_1\tilde{u}^{m-1}=(\lambda +C_1)\tilde{u}^{m-1} \quad \text{in} \quad  \Omega,\\
		\tilde{u}&=0 \quad \text{on} \quad \partial \Omega.
	\end{split} 
	\right.
\end{equation}  
By the uniqueness of the first eigenfunction of $-\Delta_m$, $\tilde u$ is the first eigenfunction.  \hfill $ \Box$

Due to Lemma $\ref{lem3}$, we consider the following equation
\begin{equation}
	\label{eq27}
	\left\{
	\begin{split}
		-\Delta_m u&=\lambda u^{m-1}+\mu u^{p-1}\quad \text{in} \quad  \Omega,\\
		u&=0 \quad \text{on} \quad \partial \Omega,
	\end{split}
	\right.
\end{equation}
where  $\mu\in[0,\hat{\mu}]$, $\hat{\mu}>0$, $p\in [m,m+\delta_0]$,  $\delta_0$  defined in Lemma $\ref{lem2}$.

The following lemma is due to Brasco and Lindgren, see \cite[Lemma 2.2]{brasco2023uniqueness}.
\begin{lemma}
	\label{lem 4}
	Let $2<m<p<m+\delta_0$ and let $\Omega\subset \mathbb{R}^N$ be an open bounded connected set. $\lambda >0$. Let $u, v\in W^{1,m}_0(\Omega)$ be two  distinct positive solutions of equation $(\ref{eq27})$. Then we must have 
	\[|\{x\in\Omega:u(x)>v(x)\}|>0 \quad \text{and} \quad |\{x\in\Omega:u(x)<v(x)\}|>0. 
	\]
	In other words, the difference $u-v$ must change sign in $\Omega$.
\end{lemma}
\noindent \emph{Proof}. We argue by contradiction and suppose, for example, that
\[|\{x\in\Omega:u(x)<v(x)\}|=0.
\]
Thus we assume that 
\[u(x)\ge v(x) \quad \text{for} \quad a.e. \quad x\in \Omega \quad \text{and} \quad |\{x\in\Omega:u(x)>v(x)\}|>0. 
\]
Let $\{v_k\}_{k\in \mathbb{N}}\subset C^{\infty}_0(\Omega)$ be a non-negative sequence  such that
\[\lim_{k\rightarrow +\infty}||\nabla v_k-\nabla v||_{L^m(\Omega)}=0.
\]
For every $\epsilon>0$, we take the admissible test function $\varphi=v^m_k/(u+\epsilon)^{m-1}$ in the weak formulation of the equation for $u$. This yields
\begin{equation}
	\begin{aligned}
		\int_{\Omega}(\lambda u^{m-1}+\mu u^{p-1})\frac{v^m_k}{(u+\epsilon)^{m-1}}dx&=\int_{\Omega}\left \langle |\nabla u|^{m-2}\nabla u,\nabla (\frac{v^m_k}{(u+\epsilon)^{m-1}})\right \rangle dx\\
		&=\int_{\Omega}\left \langle |\nabla (u+\epsilon)|^{m-2}\nabla (u+\epsilon),\nabla (\frac{v^m_k}{(u+\epsilon)^{m-1}})\right \rangle dx\\
		&\le \int_{\Omega}|\nabla v_k|^mdx.
		\nonumber
	\end{aligned}
\end{equation}
In the last inequality, we used Picone's inequality for the m-Laplacian, see \cite[Theorem 1.1]{allegretto1998picone}. By taking the limit as $k$ goes to $\infty$, we thus get 
\begin{equation}
	\int_{\Omega}(\lambda u^{m-1}+\mu u^{p-1})\frac{v^m}{(u+\epsilon)^{m-1}}dx\le \int_{\Omega}|\nabla v|^mdx=\int_{\Omega}\lambda v^m+\mu v^pdx.           \nonumber
\end{equation}
By taking the limit as $\epsilon \rightarrow 0$, using Fatou's lemma, we get
\begin{equation}
	\int_{\Omega}(\lambda +\mu u^{p-m})v^mdx\le \int_{\Omega}|\nabla v|^mdx=\int_{\Omega}\lambda v^m+\mu v^pdx.          \nonumber
\end{equation}
That is 
\begin{equation}
	\int_{\Omega}( u^{p-m}-v^{p-m})v^mdx\le 0 .    \nonumber
\end{equation}
We get a contradiction.  \hfill $\Box$
\begin{theorem}
	\label{thm5}
     Let $\delta_1>0$ and let $\Omega\subset \mathbb{R}^N$ be an open bounded connected set, with boundary of class $C^{1,\alpha}$ for some $\alpha \in (0,1)$. For every $m\in[1+\delta_1,1+\frac{1}{\delta_1}]$, and $m\le p\le m+\delta_0$, let $u_{m,p}\in W^{1,m}_0(\Omega)$ be a positive solution of $(\ref{eq27})$ with $||u_{m,p}||_{L^{\infty}(\Omega)}= 1$. Then there exist $\chi=\chi(\lambda, \hat{\mu}, \alpha, N, \delta_0,\delta_1, \Omega)\in (0,1)$, $\tau =\tau(\lambda, \hat{\mu}, \alpha, N, \delta_0, \delta_1, \Omega)>0$, and $\mu_0=\mu_0(\lambda, \hat{\mu}, \alpha, N, \delta_0, \delta_1, \Omega)>0$, $\mu_1=\mu_1(\lambda, \hat{\mu}, \alpha, N, \delta_0, \delta_1, \Omega)>0$ such that:
	
	(1) $u_{m,p}\in C^{1, \chi}(\overline{\Omega})$ with the uniform estimate
	\begin{equation}
		||u_{m,p}||_{C^{1, \chi}(\overline{\Omega})}\le L, \nonumber
	\end{equation}
    for some $L=L(\lambda, \hat{\mu}, \alpha, N, \delta_0, \delta_1, \Omega)>0$;
    
    (2) by defining $\Omega_{\tau}=\{x\in \Omega:\text{dist}(x,\p\Omega)\le \tau\}$, we have 
    \begin{equation}
    	|\nabla u_{m,p}|\ge \mu_0, \quad \text{in} \quad \Omega_{\tau}, \nonumber
    \end{equation}
    and
    \begin{equation}
    	 u_{m,p}\ge \mu_1, \quad \text{in} \quad \overline{\Omega \backslash \Omega_{\tau}}. \nonumber
    \end{equation}
\end{theorem}
\noindent \emph{Proof.} (1) is trivial. By Hopf's Lemma (see \cite[Theorem 1]{mikayelyan2015hopf}) and compactness argument. One knows
\begin{equation}
	\min_{\p \Omega} |\nabla u_{m,p}|>0, \quad \text{for  every} \quad m\in[1+\delta_1,1+\frac{1}{\delta_1}] , m\le p\le m+\delta_0. \nonumber
\end{equation}
By using that the family $\{|\nabla u_{m,p}|\}_{m\in[1+\delta_1,1+\frac{1}{\delta_1}], m\le p\le m+\delta_0}$ has a uniform $C^{0,\chi}(\p \Omega)$ estimate, an application of Arzel$\grave{a}$-Ascoli Theorem gives that there exists a constant $\overline{\mu}>0$ such that
\begin{equation}
	\min_{\p \Omega} |\nabla u_{m,p}|\ge \overline{\mu}, \quad \text{for  every} \quad m\in[1+\delta_1,1+\frac{1}{\delta_1}]  , m\le p\le m+\delta_0. \nonumber
\end{equation}
We choose $\tau_0>0$ sufficiently small, such that each point $x\in \Omega_{\tau_0}$ can be uniquely written as 
\begin{equation}
	x=x'-|x'-x|\nu_{\Omega}(x'), \quad \text{with} \quad x'\in \p \Omega. \nonumber
\end{equation}
 Here $\nu_{\Omega}$ stands for the normal outer vector. We then get for every $m\in[1+\delta_1,1+\frac{1}{\delta_1}]$, every $m\le p\le m+\delta_0$, every $0<\tau<\tau_0$ and every $x\in \Omega_{\tau}$
\begin{equation}
	|\nabla u_{m,p}(x)|\ge |\nabla u_{m,p}(x')|- ||\nabla u_{m,p}(x')|-|\nabla u_{m,p}(x)||\ge \overline{\mu}-L|x'-x|^{\chi}\ge (\overline{\mu}-L\delta^{\chi}).\nonumber
\end{equation}
If we now thoose 
\begin{equation}
	\tau =\text{min}\left\{(\frac{\overline{\mu}}{2L})^{\frac{1}{\chi}},\tau_0\right \}, \nonumber
\end{equation}
and set $\mu_0=\overline{\mu}/2$, we obtain
\begin{equation}
	|\nabla u_{m,p}(x)|\ge \mu_0 \quad \text{for every} \quad x\in \Omega_{\tau},m\in[1+\delta_1,1+\frac{1}{\delta_1}], m\le p \le m+\delta_0. \nonumber
\end{equation}
Finally, the uniform lower bound on $u_{m,p}$ in $\overline{\Omega \backslash \Omega_{\tau}}$ can be proved by observing that 
\begin{equation}
	\min_{\overline{\Omega\backslash \Omega_{\tau}}} u_{m,p}>0, \quad \text{for every} \quad m\in[1+\delta_1,1+\frac{1}{\delta_1}], m\le p \le m+\delta_0,
\end{equation}
thanks to the maximum principle. By compactness argument, one knows
\begin{equation}
	\min_{\overline{\Omega\backslash \Omega_{\tau}}} u_{m,p}\ge \mu_1, \quad \text{for every} \quad m\in[1+\delta_1,1+\frac{1}{\delta_1}], m\le p \le m+\delta_0.
\end{equation}
This conclude the proof.  \hfill $ \Box$
\begin{theorem}
	\label{thm6} (Uniform weighted Sobolev inequality). Let $\Omega\subset \mathbb{R}^N$ be an open bounded connected set, with boundary of class $C^{1,\alpha}$, for some $0<\alpha<1$. For every $m\in[1+\delta_1,1+\frac{1}{\delta_1}]$ and $m\le p\le m+\delta_0$, let $u_p\in W^{1,m}_0(\Omega)$ be a positive solution of $(\ref{eq27})$ with $||u_p||_{L^{\infty}(\Omega)}= 1$. We define 
	\begin{equation}
		\sigma_0=2\left(1-\frac{\delta_1}{2N}\right)^{-1},
	\end{equation}
then for every $2<\sigma<\sigma_0$, there exists $\mathcal{T}=\mathcal{T}(\alpha, N, \lambda, \delta_0, \delta_1, \sigma, \Omega, \hat{\mu})>0$ such that
\begin{equation}
	\label{eq34}
	\mathcal{T}\left(\int_{\Omega}|\phi|^{\sigma}dx\right)^{\frac{2}{\sigma}}\le \int_{\Omega}|\nabla u_p|^{m-2}|\nabla \phi|^2dx, \quad \text{for every} \quad \phi\in C^1(\overline{\Omega})\cap W^{1,m}_0(\Omega), \quad p\in [m,m+\delta_0].
\end{equation} 
\end{theorem}
\noindent \emph{Proof.} For every $\phi\in C^1(\overline{\Omega})\cap W^{1,m}_0(\Omega)$, we recall the classical representation formula
\begin{equation}
	\nonumber
	\phi(x)=C\int_{\mathbb{R}^N}\left\lan \nabla \phi (y),\frac{x-y}{|x-y|^N}\right\ran dy,
\end{equation}
where $C=C(N)>0$, see for example \cite[Lemma 7.14]{gilbarg1977elliptic}. This in turn implies that 
 \begin{equation}
 	\nonumber
 	|\phi(x)|\le C\int_{\mathbb{R}^N}\frac{|\nabla \phi (y)|}{|x-y|^{N-1}} dy,\quad \text{for every} \quad x\in \Omega.
 \end{equation}
We  use the notation 
\begin{equation}
	\nonumber
	t=1+\delta_2, \quad \gamma=N-1-t=N-2-\delta_2,
\end{equation}
where $\delta_2>0$, $\frac{\delta_2}{1+\delta_2}<\delta_1$.
By suitably using H\"older's inequality with exponents
\begin{equation}
	\nonumber
	2t\quad \text{and} \quad \frac{2t}{2t-1},
\end{equation}
we obtain 
\begin{equation}
	\nonumber
    |\phi(x)|\le C\left(\int_{\Omega}\frac{1}{|\nabla u_p|^{t(m-2)}|x-y|^{\gamma}}dy\right)^{\frac{1}{2t}}\left(\int_{\Omega}\left(\frac{|\nabla \phi(y)||\nabla u_p|^{\frac{m-2}{2}}}{|x-y|^{N-1-\frac{\gamma}{2t}}}\right)^{\frac{2t}{2t-1}}dy\right)^{\frac{2t-1}{2t}}.
\end{equation}
Observe that by definition 
\begin{equation}\nonumber
	t(m-2)<m-1 \quad \text{and} \quad \gamma <N-2,
\end{equation}
thus we can apply Theorem $A.5$ with $r=t(m-2)$ and get 
\begin{equation}
	\label{eq28}
	|\phi(x)|\le C\mathcal{S}^{\frac{1}{2t}}\left(\int_{\Omega}\left(\frac{|\nabla \phi(y)||\nabla u_p|^{\frac{m-2}{2}}}{|x-y|^{N-1-\frac{\gamma}{2t}}}\right)^{\frac{2t}{2t-1}}dy\right)^{\frac{2t-1}{2t}}.
\end{equation}
For simplicity, we now set 
\begin{equation}
	\nonumber
	F(y)=\left(|\nabla \phi(y)||\nabla u_p|^{\frac{m-2}{2}}\right)^{\frac{2t}{2t-1}},
\end{equation}
and observe that 
\begin{equation}
	\label{eq31}
	||F||^{\frac{2t-1}{2t}}_{L^{\frac{2t-1}{t}}(\Omega)}=\left(\int_{\Omega}|F|^{\frac{2t-1}{t}}dy\right)^{\frac{1}{2}}=\left(\int_{\Omega}|\nabla \phi|^2|\nabla u_p|^{m-2}dy\right)^{\frac{1}{2}}.
\end{equation}
We introduce the exponent $0<\Theta<N$ given by 
\begin{equation}
	\nonumber
	\left(N-1-\frac{\gamma}{2t}\right)\frac{2t}{2t-1}=N-\Theta \quad \text{that is } \quad \Theta=N-\left(N-1-\frac{\gamma}{2t}\right)\frac{2t}{2t-1}.
\end{equation}
Thanks to the choice of $\gamma$ and $t$, it is not diffcult to see that $\Theta$ is positive. More precisely, observe that this exponent is explicitly given by 
\begin{equation}
	\Theta=\frac{t-1}{2t-1}=\frac{\delta_2}{1+2\delta_2}.
\end{equation}
In view of these definitions, we can rewrite $(\ref{eq28})$ as 
\begin{equation}
	\label{eq28.5}
	|\phi(x)|\le C\mathcal{S}^{\frac{1}{2t}}\left(\int_{\Omega}\frac{|F(y)|}{|x-y|^{N-\Theta}}dy\right)^{\frac{2t-1}{2t}},
\end{equation}
so that one can recognize a suitable Riesz potential on the right hand-side. We then recall the classical potential estimate (see for example \cite[Lemma 7.12]{gilbarg1977elliptic})
\begin{equation}
	\label{eq29}
	\bigg|\bigg|\int_{\Omega}\frac{|F(y)|}{|\cdot -y|^{N-\Theta}}dy\bigg|\bigg|_{L^l(\Omega)}\le \left(\frac{1-\delta}{\frac{\Theta}{N}-\delta}\right)^{1-\delta}\omega_N^{\frac{N-\Theta}{N}}|\Omega|^{\frac{\Theta}{N}-\delta}||F||_{L^s(\Omega)},
\end{equation}
where
\begin{equation}
	\nonumber
	0\le \delta:=\frac{1}{s}-\frac{1}{l}<\frac{\Theta}{N}.
\end{equation}
By direct computation, for  $\sigma \in (2,\sigma_0)$, there holds
\begin{equation}
	\nonumber
	\frac{\sigma-2}{\sigma}\frac{t}{2t-1}<\frac{\Theta}{N}.
\end{equation}
Then,
\begin{equation}
	\nonumber
	\begin{split}
		||\phi||_{L^{\sigma}(\Omega)}&\le C\mathcal{S}^{\frac{1}{2t}}\bigg|\bigg|\left(\int_{\Omega}\frac{|F(y)|}{|\cdot-y|^{N-\Theta}}dy\right)^{\frac{2t-1}{2t}}\bigg|\bigg|_{L^{\sigma}(\Omega)}\\
		&=C\mathcal{S}^{\frac{1}{2t}}\bigg|\bigg|\int_{\Omega}\frac{|F(y)|}{|\cdot-y|^{N-\Theta}}dy\bigg|\bigg|_{L^{\sigma\frac{2t-1}{2t}(\Omega)}}^{\frac{2t-1}{2t}}.
	\end{split}
\end{equation}
Let
\begin{equation}
	\nonumber
	l=\sigma \frac{2t-1}{2t}, \quad s=\frac{2t-1}{t}.
\end{equation}
These are feasible, since
\begin{equation}
	\nonumber
	\delta=\frac{1}{s}-\frac{1}{l}=\frac{t}{2t-1}-\frac{2t}{\sigma(2t-1)}=\frac{\sigma-2}{\sigma} \frac{t}{2t-1},
\end{equation}
is positive and smaller than $\Theta/N$, thanks to the choice of $\sigma$. We then obtain
\begin{equation}
	\nonumber
	||\phi||_{L^{\sigma}(\Omega)}\le C \mathcal{S}^{\frac{1}{2t}}\left[\left(\frac{1-\delta}{\frac{\Theta}{N}-\delta}\right)^{1-\delta}\omega_N^{\frac{N-\Theta}{N}}|\Omega|^{\frac{\Theta}{N}-\delta}\right]^{\frac{2t-1}{2t}}||F||^{\frac{2t-1}{2t}}_{L^{\frac{2t-1}{t}}(\Omega)}.
\end{equation}
By recalling $(\ref{eq31})$, we thus obtained 
\begin{equation}
	\nonumber
	\sqrt{\mathcal{T}}||\phi||_{L^{\sigma}(\Omega)}\le\left(\int_{\Omega}|\nabla u_p|^{m-2}|\nabla \phi|^2dx\right)^{\frac{1}{2}},
\end{equation}
with the constant $\mathcal{T}$ given by
\begin{equation}
	\nonumber
	\mathcal{T}=\frac{1}{C^2\mathcal{S}^{\frac{1}{t}}}\left[\left(\frac{1-\delta}{\frac{\Theta}{N}-\delta}\right)^{1-\delta}\omega_N^{\frac{N-\Theta}{N}}|\Omega|^{\frac{\Theta}{N}-\delta}\right]^{\frac{1-2t}{t}}.
\end{equation}
 \hfill $\Box$
\begin{corollary}
	\label{cor1}
	Let $\delta_1\in (\frac{1}{2},1)$ and let $\Omega \subset \mathbb{R}^N$ be an open bounded connected set, with boundary of class $C^{1,\alpha}$, for some $0<\alpha <1$. For every $m\in[1+\delta_1,1+\frac{1}{\delta_1}]$ and $m\le p\le m+\delta_0$, let $u_p\in W^{1,m}_0(\Omega)$ be a positive solution of $(\ref{eq27})$ and $||u_p||_{L^{\infty}(\Omega)}=1$. Then there exists an exponent $\theta \in (1,2)$ and a constant $C=C(N, \delta_0, \delta_1, \lambda, \alpha, \Omega, \hat{\mu})>0$ such that 
	\begin{equation}
		\nonumber
		||\phi||_{W^{1,\theta}(\Omega)}\le C\left(\int_{\Omega}|\nabla u_p|^{m-2}|\nabla \phi|^2dx\right)^{\frac{1}{2}}, \quad \text{for every}\quad \phi \in W^{1,m}_0(\Omega), p\in [m,m+\delta_0].
	\end{equation}
\end{corollary}
\noindent \emph{Proof.} We take $1<\theta <(2m-2)/(2m-3)$. Then H\"older's inequality with conjugate exponents
\begin{equation}
	\nonumber
	\frac{2}{\theta} \quad \text{and} \quad \frac{2}{2-\theta},
\end{equation}
imply
\begin{equation}
	\nonumber
	\begin{split}
		\left(\int_{\Omega}|\nabla \phi|^\theta dx\right)^{\frac{1}{\theta}}&\le \left(\int_{\Omega}|\nabla u_p|^{m-2}|\nabla \phi|^2dx\right)^{\frac{1}{2}}\left(\int_{\Omega}\frac{1}{|\nabla u_p|^{\frac{\theta}{2-\theta}(m-2)}}dx\right)^{\frac{2-\theta}{2\theta}}\\
		&\le \tilde{\mathcal{S}}^{\frac{2-\theta}{2\theta}}\left(\int_{\Omega}|\nabla u_p|^{m-2}|\nabla \phi|^2dx\right)^{\frac{1}{2}},
	\end{split}
\end{equation}
where we have used that
\begin{equation}
	\nonumber
	\frac{\theta}{2-\theta}(m-2)<(m-1),
\end{equation}
which allows us to use estimate $(A.14)$ from Theorem A.5, with $r=(m-2)\theta/(2-\theta)$. Since $\Omega$ is bounded, by H\"older's inequality we have that $W^{1,m}_0(\Omega)\subset W^{1,\theta}_0(\Omega)$, with continuous inclusion. Moreover, by Poincar\'e inequality
\begin{equation}
	\nonumber
	||\nabla \phi||_{L^{\theta}(\Omega)} \quad \text{and} \quad ||\phi||_{W^{1,\theta}(\Omega)},
\end{equation}
are equivalent norms on $W^{1,\theta}_0(\Omega)$. These facts conclude the proof. \hfill $\Box$
\begin{corollary}
	\label{cor2}
	(Uniform compact embedding) Let $m\in[1+\delta_1,1+\frac{1}{\delta_1}]$ and let $\Omega \subset \mathbb{R}^N$ be an open bounded connected set, with boundary of class $C^{1,\alpha}$, for some $0<\alpha <1$. We take a sequence $\{p_k\}_{{}k\in \mathbb{N}}\subset [m,m+\delta_0]$ and consider accordingly $u_k\in W^{1,m}_0(\Omega)$ a positive solutions of $(\ref{eq27})$ with $p=p_k$ and $||u_k||_{L^{\infty}(\Omega)}=1$. If $\{\phi _k\}_{k\in \mathbb{N}}\subset W^{1,m}_0(\Omega)\cap C^1(\overline{\Omega})$ is a sequence of functions satisfying
	\[\int_{\Omega}|\nabla u_k|^{m-2}|\nabla \phi_k|^2dx\le C,\quad \text{for  every} \quad k \in \mathbb{N},
	\]
	then $\{\phi_k\}_{k\in \mathbb{N}}$ converges strongly in $L^2(\Omega)$ and weakly in $W^{1,\theta}_0(\Omega)$, up to a subsequence. Here $\theta$ is the same exponent as in Corollary $\ref{cor1}$. 
\end{corollary}
\noindent \emph{Proof.} The assumption, in conjunction with Corollary $\ref{cor1}$, entails that $\{\phi_k\}_{k\in\mathbb{N}}$ is a bounded sequence in $W^{1,\theta}_0(\Omega)$. By the classical Rellich-Kondra\v{s}ov Theorem, we get that this sequence converges weakly in $W^{1,\theta}(\Omega)$ and strongly in $L^{\theta}(\Omega)$, up to a subsequence. Moreover, since $W^{1,\theta}_0(\Omega)$ is also weakly closed, we get that the limit still belongs to $W^{1,\theta}_0(\Omega)$.

In order to get the strong $L^2$ convergence, we observe that, if we denote by $\sigma_0$ the exponent of Theorem $\ref{thm6}$, for every $2<\sigma <\sigma_0$ and $n,k \in \mathbb{N}$ we have
\begin{equation}
	\nonumber
		||\phi_n-\phi_k||_{L^2(\Omega)}\le 	||\phi_n-\phi_k||^{1-\tau}_{L^{\sigma}(\Omega)} ||\phi_n-\phi_k||^{\tau}_{L^{\theta}(\Omega)}
		\le C ||\phi_n-\phi_k||^{\tau}_{L^{\theta}(\Omega)}.
\end{equation} 
This ends the proof of the corollary.    \hfill $\Box$

Similarly, we have the following corollary.
\begin{corollary}
	\label{cor2.1}
	(Uniform compact embedding) Let $\Omega \subset \mathbb{R}^N$ be an open bounded connected set, with boundary of class $C^{1,\alpha}$, for some $0<\alpha <1$. We take a sequence $\{m_k\}_{{}k\in \mathbb{N}}\subset [2,2+\delta_1]$ and consider accordingly $u_k\in W^{1,m_k}_0(\Omega)$ a positive solutions of $(\ref{eq27})$ with $m=m_k$ and $||u_k||_{L^{\infty}(\Omega)}=1$. If $\{\phi _k\}_{k\in \mathbb{N}}\subset W^{1,m_k}_0(\Omega)\cap C^1(\overline{\Omega})$ is a sequence of functions satisfying
	\[\int_{\Omega}|\nabla u_k|^{m_k-2}|\nabla \phi_k|^2dx\le C,\quad \text{for  every} \quad k \in \mathbb{N},
	\]
	then $\{\phi_k\}_{k\in \mathbb{N}}$ converges strongly in $L^2(\Omega)$ and weakly in $W^{1,\theta}_0(\Omega)$, up to a subsequence. Here $\theta$ is the same exponent as in Corollary $\ref{cor1}$. 
\end{corollary}
\section{ Eigenvalue problem}
In this section, we treat a weighted linear eigenvalue problem.
It is convenient to introduce the notation 
\begin{equation}
	\nonumber
	\mathcal{H}(z)=\frac{1}{m}|z|^m, \quad \text{for every} \quad z\in \mathbb{R}^N.
\end{equation}
Then we observe that
\begin{equation}
	\nonumber
	\nabla \mathcal{H}(z)=|z|^{m-2}z, \quad \text{for every} \quad z\in \mathbb{R}^N,
\end{equation}
and 
\begin{equation}
	\label{ep1}
	D^2 \mathcal{H}(z)=|z|^{m-2}\text{Id}+(m-2)|z|^{m-4}z\otimes z, \quad \text{for every} \quad z\in \mathbb{R}^N.
\end{equation}
In particular, we have the following facts 
\begin{equation}
	\label{ep2}
	D^2\mathcal{H}(z)z=(m-1)|z|^{m-2}z, \quad \lan D^2\mathcal{H}(z)\xi,\xi \ran\le (m-1)|z|^{m-2}|\xi|^2, \quad \text{for} \quad z,\xi\in \mathbb{R}^N.
\end{equation}
We have the following elementary inequality.
\begin{lemma}
	\label{lem4}
	Let $2<m<+\infty$ and let $\Omega\subset \mathbb{R}^N$ be an open set. For every $v, w, \phi\in W^{1,1}_{\mathrm{loc}}(\Omega)$, we have
	\begin{equation}
		\nonumber
		\begin{split}
			&|\lan D^2\mathcal{H}(\nabla \phi)\nabla v,\nabla v\ran-\lan D^2\mathcal{H}(\nabla \phi)\nabla w,\nabla w\ran|\\
			&\le (m-1)|\nabla \phi|^{m-2}|\nabla v-\nabla w|(|\nabla v|+|\nabla w|),\quad a.e. \quad \text{on} \quad \Omega.
		\end{split}
	\end{equation}
\end{lemma}
\noindent \emph{Proof.} By Lemma A.2 of \cite{brasco2023uniqueness}, we have 
\begin{equation}
	\nonumber
	|\lan D^2\mathcal{H}(\nabla \phi)\nabla v,\nabla v\ran-\lan D^2\mathcal{H}(\nabla \phi)\nabla w,\nabla w\ran|\le|D^2\mathcal{H}(\nabla \phi)(\nabla v-\nabla w)|(|\nabla v|+|\nabla w|).
\end{equation}
Since the Hessian matrix is given by $(\ref{ep1})$, we get
\begin{equation}
	\nonumber
	|D^2\mathcal{H}(\nabla \phi)(\nabla v-\nabla w)|\le (m-1)|\nabla \phi|^{m-2}|\nabla v-\nabla w|.
\end{equation}
By using this inequality in the first estimate, we conclude the proof.\hfill $\Box$
\begin{proposition}
	\label{pro3.2}
	Let $2<m<+\infty$ and let $\Omega \subset \mathbb{R}^N$ be an open connected set with finite volume. Let $U\in W^{1,m}_0(\Omega)$ be the unique positive extremal of
	\begin{equation}
		\label{ep3}
		\lambda_m(\Omega)=\inf_{\phi\in W^{1,m}_0(\Omega)}\left\{\int_{\Omega}|\nabla \phi|^mdx:\int_{\Omega}|\phi|^mdx=1\right\}.
	\end{equation}
   By setting
   \begin{equation}
   	\nonumber
   	\lambda(\Omega;U)=\inf_{\phi\in C^\infty_0(\Omega)}\left\{\int_{\Omega}\lan D^2\mathcal{H}(\nabla U)\nabla \phi,\nabla \phi\ran dx:\int_{\Omega}U^{m-2}|\phi|^2dx=1\right\},
   \end{equation}
    we have 
    \begin{equation}
    	\nonumber
    	\lambda(\Omega;U)=(m-1)\lambda_m(\Omega).
    \end{equation}
\end{proposition}
\noindent \emph{Proof.} The inequality 
\begin{equation}
	\label{ep4}
	\lambda(\Omega;U)\le (m-1)\lambda_m(\Omega),
\end{equation}
is straightforward. Indeed, for every $\epsilon>0$, we take a non-negative $U_{\epsilon}\in C^\infty_0(\Omega)$ such that
\begin{equation}
	\nonumber
	\int_{\Omega}|\nabla U_{\epsilon}|^mdx\le \lambda_{m}(\Omega)+\epsilon, \int_{\Omega}U^m_{\epsilon}dx=1 \quad \text{and} \quad \lim_{\epsilon\rightarrow 0}\int_{\Omega}|\nabla U_{\epsilon}-\nabla U|^mdx=0.
\end{equation}
Thus in particular we have convergence in $L^m(\Omega)$, as well. Denote 
\begin{equation}
	\nonumber
	\tilde{U}_\epsilon :=\frac{U_\epsilon}{(\int_{\Omega}U^{m-2}U_{\epsilon}^2dx)^{\frac{1}{2}}}.
\end{equation}
It is easy to see that
\begin{equation}
	\label{ep5}
	\lim_{\epsilon \rightarrow 0}\int_{\Omega}U^{m-2}|\tilde U_{\epsilon}|^2dx=\int_{\Omega}U^mdx=1.
\end{equation}
Also
\begin{equation}
	\label{ep6}
	\lim_{\epsilon\rightarrow 0}\int_{\Omega}\lan D^2\mathcal{H}(\nabla U)\nabla U_\epsilon,\nabla U_\epsilon \ran dx=\int_{\Omega}\lan D^2\mathcal{H}(\nabla U)\nabla U,\nabla U \ran dx.
\end{equation}
Indeed, by applying Lemma $\ref{lem4}$, we get 
\begin{equation}
	\nonumber
	\begin{split}
		\bigg|\int_{\Omega}\lan D^2\mathcal{H}(\nabla U)\nabla U_\epsilon,\nabla U_\epsilon \ran dx-&\int_{\Omega}\lan D^2\mathcal{H}(\nabla U)\nabla U,\nabla U \ran dx\bigg|\\
		\le&(m-1)\int_{\Omega}|\nabla U|^{m-2}|\nabla U_\epsilon-\nabla U|(|\nabla U_\epsilon|+|\nabla U|)dx.
	\end{split}
\end{equation}
Then $(\ref{ep6})$ follows by using H\"older's inequality.
By $(\ref{ep2})$, $(\ref{ep5})$ and $(\ref{ep6})$,  
we get 
\begin{equation}
	\nonumber
	\lambda(\Omega;U)\le (m-1)\lambda_m(\Omega), \quad \text{as} \quad \epsilon \rightarrow 0.
\end{equation}
For the converse inequality, we first recall that $u$ satisfies 
\begin{equation}
	\label{ep7}
	\int_{\Omega}\lan|\nabla U|^{m-2}\nabla U, \nabla \phi\ran dx=\lambda_m(\Omega)\int_{\Omega}U^{m-1}\phi dx, \quad \text{for every} \quad \phi\in W^{1,m}_0(\Omega), 
\end{equation}
by minimality. By using $(\ref{ep2})$, this can be also written as 
\begin{equation}
	\label{ep8}
	\int_{\Omega}\lan\D^2\mathcal{H}(\nabla U)\nabla U, \nabla \phi\ran dx=(m-1)\lambda_m(\Omega)\int_{\Omega}U^{m-1}\phi dx, \quad \text{for every} \quad \phi\in W^{1,m}_0(\Omega). 
\end{equation}
We take $\epsilon>0$ and $\phi_\epsilon\in C^\infty_0(\Omega)$ such that
\begin{equation}
	\nonumber
	\int_{\Omega}\lan\D^2\mathcal{H}(\nabla U)\nabla \phi_\epsilon, \nabla \phi_\epsilon\ran dx<\lambda(\Omega;U)+\epsilon \quad \text{and} \quad \int_{\Omega}U^{m-2}\phi_\epsilon^2dx=1.
\end{equation}
Then we insert the test function $\phi_\epsilon ^2/U$ in the equation $(\ref{ep8})$, so to get
\begin{equation}
	\label{ep9}
	(m-1)\lambda_m(\Omega)\int_{\Omega}U^{m-1}\frac{\phi_\epsilon^2}{U} dx=\int_{\Omega}\left\lan\D^2\mathcal{H}(\nabla U)\nabla U, \nabla (\frac{\phi_\epsilon^2}{U})  \right\ran dx.
\end{equation}
We now use Picone's identity of Lemma A.1  in \cite{brasco2023uniqueness} with the choice $A=D^2\mathcal{H}(\nabla U)$. This gives 
\begin{equation}
	\nonumber
	\begin{split}
		\left \lan\D^2\mathcal{H}(\nabla U)\nabla U, \nabla (\frac{\phi_\epsilon^2}{U}) \right \ran&=\left \lan \D^2\mathcal{H}(\nabla U)\nabla \phi_\epsilon, \nabla \phi_\epsilon \right\ran\\
		&-\left \lan \D^2\mathcal{H}(\nabla U)\left(\phi_\epsilon\frac{\nabla U}{U}-\nabla \phi_\epsilon\right),\left(\phi_\epsilon\frac{\nabla U}{U}-\nabla \phi_\epsilon\right) \right \ran.
	\end{split}
\end{equation}
By  $(\ref{ep9})$, we get 
\begin{equation}
	\nonumber
	\begin{split}
		(m-1)\lambda_m(\Omega)\int_{\Omega}U^{m-2}\phi_\epsilon^2 dx&=\int_{\Omega}\lan \D^2\mathcal{H}(\nabla U)\nabla \phi_\epsilon, \nabla \phi_\epsilon \ran dx\\
		&-\int_{\Omega}\left\lan \D^2\mathcal{H}(\nabla U)\left(\phi_\epsilon\frac{\nabla U}{U}-\nabla \phi_\epsilon\right ),\left(\phi_\epsilon\frac{\nabla U}{U}-\nabla \phi_\epsilon\right) \right\ran dx\\
		&\le\lambda (\Omega,U)+\epsilon,
	\end{split}
\end{equation}
thanks to the fact that $D^2\mathcal{H}(\nabla U)$ is positive semidefinite. By recalling that 
\begin{equation}
	\nonumber
	\int_{\Omega}U^{m-2}\phi_\epsilon ^2dx=1,
\end{equation}
and using the arbitrariness of $\epsilon >0$, we finally get the desired conclusion. \hfill $\Box$
\begin{definition}
	For $m>2$, we define the weighted Sobolev space
	\begin{equation}
		\nonumber
		X^{1,2}(\Omega;|\nabla U|^{m-2}):=\left\{\phi \in W^{1,1}_{\mathrm{loc}}(\Omega)\cap L^2(\Omega):\int_{\Omega}|\nabla U|^{m-2}|\nabla \phi|^2dx<+\infty\right \},
	\end{equation}
endowed with the natural norm
\begin{equation}
	\nonumber
	||\phi||_{X^{1,2}(\Omega;|\nabla U|^{m-2})}=||\phi||_{L^2(\Omega)}+\left(\int_{\Omega}|\nabla U|^{m-2}|\nabla \phi|^2dx\right)^{\frac{1}{2}}.
\end{equation}
\end{definition}
We set $X^{1,2}_0(\Omega;|\nabla U|^{m-2})$ for the completion of $C^\infty_0(\Omega)$ with respect to this norm.
\begin{lemma}
	Let $m>2$ and let $\Omega\subset \mathbb{R}^N$ be an open bounded connected set, with $C^{1,\alpha}$ boundary, for some $0<\alpha<1$. With the notation above, we have
	\begin{equation}
		\nonumber
	 X^{1,2}_0(\Omega;|\nabla U|^{m-2})\subset W^{1,1}_0(\Omega),
	\end{equation}
    with continuous inclusions.
\end{lemma}
\noindent \emph{Proof.}  It is sufficient to prove that there exists two constants $C_1>0$ such that
\begin{equation}
	\nonumber
	C_1||\phi||_{W^{1,1}(\Omega)}\le ||\phi||_{X^{1,2}(\Omega;|\nabla U|^{m-2})}, \quad \text{for every} \quad \phi\in C^\infty_0(\Omega).
\end{equation}
The estimate follows from Corollary $\ref{cor1}$, with $m=p$.  \hfill $ \Box$
\begin{remark}
	\label{re4.1}
	It is easy to verify that $U\in X^{1,2}_0(\Omega;|\nabla U|^{m-2})$.
\end{remark}
\begin{proposition}
	\label{prop3.5}
	Let $2<m<+\infty$ and let $\Omega\subset \mathbb{R}^N$ be an open bounded connected set, with $C^{1,\alpha}$ boundary, for some $0<\alpha<1$. With the notation above, the infimum $\lambda(\Omega;U)$ is uniquely attained on the space  $X^{1,2}_0(\Omega;|\nabla U|^{m-2})$ by the function $U$ or $-U$.
\end{proposition}
\noindent \emph{Proof.} We first notice that if $v\in X^{1,2}_0(\Omega;|\nabla U|^{m-2})$ and $\{v_k\}_{k\in \mathbb{N}} \subset C^\infty_0(\Omega)$ is such that
\begin{equation}
	\nonumber
	\lim_{k\rightarrow \infty}||v_k-v||_{X^{1,2}(\Omega;|\nabla U|^{m-2})}=0,
\end{equation}
then
\begin{equation}
	\label{ep10}
	\lim_{k\rightarrow \infty}\int_{\Omega}\lan D^2\mathcal{H}(\nabla U)\nabla v_k,\nabla v_k\ran dx=\int_{\Omega}\lan D^2\mathcal{H}(\nabla U)\nabla v,\nabla v\ran dx.
\end{equation}
Indeed, by Lemma $\ref{lem4}$
\begin{equation}
	\nonumber
	|\lan D^2\mathcal{H}(\nabla U)\nabla v_k,\nabla v_k\ran-\lan D^2\mathcal{H}(\nabla U)\nabla v,\nabla v\ran|\le (m-1)|\nabla U|^{m-2}|\nabla v_k-\nabla v|(|\nabla v_k|+|\nabla v|).
\end{equation}
By integrating over $\Omega$ and using H\"older's inequality, we have 
\begin{equation}
	\nonumber
	\begin{split}
		|\int_{\Omega}\lan D^2\mathcal{H}(\nabla U)\nabla v_k,\nabla v_k\ran dx&-\int_{\Omega}\lan D^2\mathcal{H}(\nabla U)\nabla v,\nabla v\ran dx|\\
		&\le C\int_{\Omega}|\nabla U|^{m-2}|\nabla v_k-\nabla v|(|\nabla v_k|+|\nabla v|)dx\\
		&\le C\left(\int_\Omega |\nabla U|^{m-2}|\nabla v_k-\nabla v|^2dx\right)^{\frac{1}{2}}\\
		&\times \left(\int_\Omega |\nabla U|^{m-2}(|\nabla v_k|+|\nabla v|)^2dx\right)^{\frac{1}{2}}.
	\end{split}
\end{equation}
By observing that the last term converges to $0$, we get $(\ref{ep10})$.
Since $C^\infty_0(\Omega)$ is dense in $X^{1,2}_0(\Omega;|\nabla U|^{m-2})$ by definition, our previous computations show that
\begin{equation}
	\nonumber
	\lambda(\Omega;U)=\inf_{\phi\in X^{1,2}_0(\Omega;|\nabla U|^{m-2})}\left\{\int_\Omega \lan D^2\mathcal{H}(\nabla U)\nabla \phi,\nabla \phi\ran:\int_\Omega U^{m-2}\phi^2dx=1 \right\}.
\end{equation}
In order to prove that any minimizer must coincide either with $U$ or $-U$, we assume that there is another minimizer $v\in X^{1,2}_0(\Omega;|\nabla U|^{m-2})$. By definition, there exists a sequence $\{v_k\}_{k\in \mathbb{N}}\subset C^\infty_0(\Omega)$ such that
\begin{equation}
	\nonumber
	\lim_{k\rightarrow +\infty}\left[\int_\Omega |\nabla U|^{m-2}|\nabla v_k-\nabla v|^2dx+\int_\Omega|v_k-v|^2dx\right]=0.
\end{equation}
We recall that $u$ satisfies $(\ref{ep8})$. For every $k\in \mathbb{N}$, the choice $\phi=v_k^2/U$ is feasible in $(\ref{ep8})$ and it yields 
\begin{equation}
	\label{ep11}
	\begin{split}
		\lambda(\Omega;U)\int_\Omega U^{m-2}v_k^2dx&=\int_\Omega \left\lan D^2\mathcal{H}(\nabla U)\nabla U,\nabla \left(\frac{v_k^2}{U}\right) \right\ran dx\\
		&= \int_\Omega \lan D^2\mathcal{H}(\nabla U)\nabla v_k,\nabla v_k \ran dx\\
		&-\int_\Omega \left\lan D^2\mathcal{H}(\nabla U)\left(v_k\frac{\nabla U}{U}-\nabla v_k\right),\left(v_k\frac{\nabla U}{U}-\nabla v_k\right)\right\ran dx,
	\end{split}
\end{equation}
where in the second equality we used the general version of Picone identity given by Lemma $A.1$ in \cite{brasco2023uniqueness}, with the positive semi-definite matrix $A=D^2\mathcal{H}(\nabla U)$.
We now wish to pass the limit as $k$ goes to $\infty$ in the previous identity. We notice at first that 
\begin{equation}
	\nonumber
	\lim_{k\rightarrow +\infty}\int_\Omega U^{m-2}v^2_kdx=\int_\Omega U^{m-2}v^2dx=1.
\end{equation}
As for the first term on the right-hand side of $(\ref{ep11})$, we simply use $(\ref{ep10})$.
We are left with handing the last term in $(\ref{ep11})$. We have that $\{(v_k,\nabla v_k)\}_{k\in \mathbb{N}}$ converges almost everywhere to $(v,\nabla v)$, possibly up to extracting a subsequence. By observing that $\D^2\mathcal{H}(\nabla U)$ is positive semidefinite, an application of Fatou's Lemma yields
\begin{equation}
	\nonumber
	\begin{split}
		\liminf_{k\rightarrow \infty}&\int_\Omega\left\lan \D^2\mathcal{H}(\nabla U)\left(v_k\frac{\nabla U}{U}-\nabla v_k\right),\left(v\frac{\nabla U}{U}-\nabla v_k\right)\right\ran dx\\
		&\ge \int_\Omega\left\lan \D^2\mathcal{H}(\nabla U)\left(v\frac{\nabla U}{U}-\nabla v\right),\left(v\frac{\nabla U}{U}-\nabla v\right)\right\ran dx.
	\end{split}
\end{equation}
Thus, by taking the limit as $k$ goes to $\infty$ in $(\ref{ep11})$, we get
\begin{equation}
	\nonumber
	\lambda(\Omega;U)+\int_\Omega \left\lan \D^2\mathcal{H}(\nabla U)\left(v\frac{\nabla U}{U}-\nabla v\right),\left(v\frac{\nabla U}{U}-\nabla v\right)\right\ran dx \le \lambda(\Omega;U).
\end{equation}
This entails that we must have
\begin{equation}
	\label{ep12}
	\left\lan \D^2\mathcal{H}(\nabla U)\left(v\frac{\nabla U}{U}-\nabla v\right),\left(v\frac{\nabla U}{U}-\nabla v\right)\right\ran=0, \quad \text{a.e. in} \quad \Omega.
\end{equation}
From the definition of $D^2\mathcal{H}$, it is clear that $D^2\mathcal{H}(\nabla U)$ is positive definite whenever $\nabla U$ does not vanish, which is true almost everywhere by $(A.14)$. Therefore, from $(\ref{ep12})$ we must have 
\begin{equation}
	\label{ep13}
	v\frac{\nabla U}{U}-\nabla v=0, \quad \text{a.e. in} \quad \Omega.
\end{equation}
We now observe that $v\in W^{1,1}_0(\Omega)$ thanks to Lemma $3.4$. $U \in C^{1}(\overline{\Omega})$ as in \cite[Proposition 3.5]{brasco2023uniqueness} and it has the following property: for every $\Omega'\Subset\Omega$, there exists a constant $C=C(\Omega')>0$ such that $U\ge 1/C$ on $\Omega'$. Thus we have
\begin{equation}
	\nonumber
	\frac{v}{U}\in W^{1,1}_{\mathrm{loc}}(\Omega),
\end{equation}
and Leibnitz's rule holds for its distributional gradient. This is given by
\begin{equation}
	\nonumber
	\nabla \left(\frac{v}{U}\right)=\frac{U\nabla v-v\nabla U}{U^2}, \quad \text{a.e. in} \quad \Omega,
\end{equation}
and thus it identically vanishes almost everywhere in $\Omega$, by virtue of $(\ref{ep13})$. Since $\Omega$ is connected, this implies that $v/U$ is constant in $\Omega$. Thus we get that $v$ is proportional to $U$ in $\Omega$. This desired result is now a consequence of the normalization taken.  \hfill $ \Box$
\section{Proof of the Theorems 1.2}
\noindent \emph{Proof of Theorem 1.2}. We divide the proof into three parts, for ease of readability. 

\noindent \emph{Par 1}. Linearized equation. We argue by contradiction: we suppose that for $\{p_k\}_{k\in \mathbb{N}}$ with $p_k\searrow m$, the problem $(\ref{eq2})$ admits (at least) two distinct positive solutions $u_k$ and $v_k$ for $p=p_k$. 

Without loss of generality, we can assume $||u_k||_{L^\infty(\Omega)}\ge ||v_k||_{L^\infty(\Omega)}$ for each $k$. Let $M_k:=||u_k||_{L^\infty(\Omega)}$.
We may choose a subsequence such that $\lim_{k\rightarrow \infty}\frac{||v_k||_{L^\infty(\Omega)}}{M_k}=\mu'\in [0,1]$.
Define $\tilde{u}_k=\frac{u_k}{M_k}$, $\tilde{v}_k=\frac{v_k}{M_k}$. Observe that they solve 
\begin{equation}
	\nonumber
	-\Delta_m \tilde{u}_k=\lambda \tilde{u}_k^{m-1}+M^{p_k-m}_k\tilde{u}^{p_k-1}_k \quad \text{and} \quad -\Delta_m \tilde{v}_k=\lambda \tilde{v}_k^{m-1}+M^{p_k-m}_k\tilde{v}^{p_k-1}_k, \quad \text{in} \quad \Omega,
\end{equation}
with 
\[\tilde{u}_k=\tilde{v}_k=0 \quad \text{on} \quad \partial \Omega. 
\]
According to Lemma $\ref{lem3}$, we can choose a subsequence such that $\tilde{u}_k$ converges to $\tilde{u}$ in $C^{1,\beta}$, where $\tilde{u}$ is the first eigenfunction of $-\Delta_m$ in $\Omega$. In a similar way, $\tilde{v}_k$ converges to $\mu'\tilde{u}$ in $C^{1,\beta}$.
The equations solved by $\tilde{u}_k$ and $\tilde{v}_k$ can be written as
\begin{equation}
	\label{ee1}
	\int_\Omega \lan \nabla \mathcal{H}(\nabla \tilde{u}_k),\nabla \phi \ran dx=\lambda \int_\Omega \tilde{u}_k^{m-1}\phi dx+M^{p_k-m}_k\int_\Omega \tilde{u}^{p_k-1}_k\phi dx,
\end{equation}
and
\begin{equation}
	\label{ee2}
	\int_\Omega \lan \nabla \mathcal{H}(\nabla \tilde{v}_k),\nabla \phi \ran dx=\lambda \int_\Omega \tilde{v}_k^{m-1}\phi dx+M^{p_k-m}_k\int_\Omega \tilde{v}^{p_k-1}_k\phi dx,
\end{equation}
for any $\phi \in W^{1,m}_0(\Omega)$.
 
We now observe that for every $z,w\in \mathbb{R}^N$ we have 
\begin{equation}
	\label{ee3}
	\begin{split}
		\nabla \mathcal{H}(z)-\nabla \mathcal{H}(w)&=\int_{0}^{1}\frac{d}{dt}(\nabla \mathcal{H}(tz+(1-t)w)dt)\\
		&=\left \lan \int_{0}^{1}(D^2 \mathcal{H}(tz+(1-t)w))dt,z-w \right \ran.
	\end{split}
\end{equation}
Similarly, for every $a,b\in[0,+\infty)$ we have
\begin{equation}
	\label{ee4}
	\begin{split}
		a^p-b^p&=\int_{0}^{1}\frac{d}{dt}(ta+(1-t)b)^pdt\\
		&=p\left(\int_{0}^{1}(ta+(1-t)b)^{p-1}dt\right)(a-b).
	\end{split}
\end{equation}
By substracting the two equations $(\ref{ee1})$ and $(\ref{ee2})$, using $(\ref{ee3})$ with $z=\nabla \tilde{u}_k(x)$, $w=\nabla \tilde{v}_k(x)$
and $(\ref{ee4})$, we get 
\begin{equation}
	\label{ee5}
	\begin{split}
		\int_\Omega\lan A_k(x)\nabla (\tilde u_k-\tilde{v}_k),\nabla \phi \ran dx&=(m-1)\lambda \int_\Omega\int_{0}^{1}(t\tilde{u}_k+(1-t)\tilde{v}_k)^{m-2}(\tilde{u}_k-\tilde{v}_k)\phi dtdx\\
		&+M_k^{p_k-m}(p_k-1)\int_\Omega\int_{0}^{1}(t\tilde{u}_k+(1-t)\tilde{v}_k)^{p_k-2}(\tilde{u}_k-\tilde{v}_k)\phi dtdx,
	\end{split}
\end{equation}
where
\begin{equation}
	\nonumber
	A_k(x)=\int_{0}^{1}D^2\mathcal{H}(t\nabla \tilde{u}_k+(1-t)\nabla \tilde{v}_k)dt.
\end{equation}
Set 
\begin{equation}
	\nonumber
	w_k=\frac{\tilde u_k-\tilde{v}_k}{||\tilde u_k-\tilde{v}_k||_{L^2(\Omega)}}\in W^{1,m}_0(\Omega).
\end{equation}
Then by $(\ref{ee5})$ we get that $w_k$ satisfies 
\begin{equation}
	\label{ee6}
	\begin{split}
		\int_\Omega\lan A_k(x)\nabla w_k,\nabla \phi \ran dx&=(m-1)\lambda \int_\Omega\int_{0}^{1}(t\tilde{u}_k+(1-t)\tilde{v}_k)^{m-2}w_k\phi dtdx\\
		&+M_k^{p_k-m}(p_k-1)\int_\Omega\int_{0}^{1}(t\tilde{u}_k+(1-t)\tilde{v}_k)^{p_k-2}w_k\phi dtdx.
	\end{split}
\end{equation}
Let $\phi=w_k$. Then 
\begin{equation}
	\label{ee7}
	\begin{split}
		\int_\Omega\lan A_k(x)\nabla w_k,\nabla w_k \ran dx&=(m-1)\lambda \int_\Omega\int_{0}^{1}(t\tilde{u}_k+(1-t)\tilde{v}_k)^{m-2}w_k^2 dtdx\\
		&+M_k^{p_k-m}(p_k-1)\int_\Omega\int_{0}^{1}(t\tilde{u}_k+(1-t)\tilde{v}_k)^{p_k-2}w_k^2 dtdx.
	\end{split}
\end{equation}
\emph{Par 2.} Convergence of $w_k$. Since 
\begin{equation}
	\label{ee8}
	||w_k||_{L^2(\Omega)}=1, 
\end{equation}
one has
\begin{eqnarray}
	\nonumber
	\int_\Omega\lan A_k(x)\nabla w_k,\nabla w_k \ran dx\le C.
\end{eqnarray}
An application of Lemma $A.3$ in \cite{brasco2023uniqueness} yields
\begin{equation}
	\label{ee9}
	\begin{split}
		\lan A_k\xi,\xi\ran &=\int_{0}^{1} \left\lan D^2\mathcal{H}(t\nabla \tilde{u}_k+(1-t)\nabla \tilde{v}_k)\xi,\xi \right\ran dt\\
		&\ge\left(\int_{0}^{1}|t\nabla \tilde{u}_k+(1-t)\nabla \tilde{v}_k|^{m-2}dt\right)|\xi|^2\ge \frac{1}{4^{m-1}}(|\nabla \tilde{u}_k|+|\nabla \tilde{v}_k|)^{m-2}|\xi|^2
.	\end{split}
\end{equation}
Therefore, we obtain 
\begin{equation}
	\nonumber
	\int_\Omega (|\nabla \tilde{u}_k|+|\nabla \tilde{v}_k|)^{m-2}|\nabla w_k|^2dx\le C ,\quad \text{for every} \quad k\in\mathbb{N}.
\end{equation}
By Corollary $\ref{cor2}$, there exist $\theta=\theta(m)\in (1,2)$ and $w\in L^2(\Omega)\cap W^{1,\theta}_0(\Omega)$ such that $\{w_k\}_{k\in \mathbb{N}}$ converges strongly in $L^2(\Omega)$ and weakly in $W^{1,\theta}(\Omega)$ to $w$, up to a subsequence. In addition, from the $C^1(\overline{\Omega})$ convergence of $\tilde{u}_k$ and $\tilde{v}_k$, we have 
\begin{equation}
	\label{ee9.1}
	A_k\rightarrow \int_{0}^{1}D^2\mathcal{H}(t\nabla \tilde{u}+(1-t)\mu'\nabla \tilde{u})dt   \quad \text{uniformly on }  \Omega .
\end{equation}
This is now enough to pass the limit in $(\ref{ee6})$ and $(\ref{ee7})$. Indeed, the convergence of the right-hand side of $(\ref{ee6})$ follows from the uniform convergences of $\tilde{u}_k$ and $\tilde{v}_k$. As for the left hand side, we have for every $\phi \in C^\infty_0(\Omega)$
\begin{equation}
	\nonumber
	\begin{split}
		\bigg|\int_\Omega \lan A_k(x)\nabla w_k&,\nabla \phi\ran dx-\int_\Omega \int_{0}^{1} \lan D^2\mathcal{H}(t\nabla \tilde{u}+(1-t)\mu'\nabla \tilde{u})\nabla w,\nabla \phi\ran dx\bigg|\\
		&\le \bigg|\int_{\Omega}\left\lan \left(A_k(x)-\int_{0}^{1}\lan D^2\mathcal{H}(t\nabla \tilde{u}+(1-t)\mu'\nabla \tilde{u})dt\right)\nabla w_k,\nabla \phi \right \ran dx\bigg|\\
		&+\bigg|\int_{\Omega}\left\lan \int_{0}^{1} D^2\mathcal{H}(t\nabla \tilde{u}+(1-t)\mu'\nabla \tilde{u})dt(\nabla w_k-\nabla w),\nabla \phi\right \ran dx\bigg|\\
		&\le \bigg|\bigg|A_k- \int_{0}^{1}\lan D^2\mathcal{H}(t\nabla \tilde{u}+(1-t)\mu'\nabla \tilde{u})dt\bigg|\bigg|_{L^\infty(\Omega)}||\nabla \phi||_{L^\infty(\Omega)}\int_\Omega |\nabla w_k|dx\\
		&+\bigg|\int_{\Omega}\left\lan \int_{0}^{1} D^2\mathcal{H}(t\nabla \tilde{u}+(1-t)\mu'\nabla \tilde{u})dt(\nabla w_k-\nabla w),\nabla \phi\right\ran dx\bigg|\rightarrow 0.
	\end{split}
\end{equation}
The last estimate of the above comes from $(\ref{ee9.1})$, $||w_k||_{W^{1,\theta}(\Omega)}\le C$ and $D^2\mathcal{H}(\nabla \tilde{u})\in L^{\infty}(\Omega)$.

Thus $w$ satisfies
\begin{equation}
	\nonumber
	\begin{split}
		\int_{\Omega} \left\lan \int_{0}^{1}D^2\mathcal{H}(t\nabla \tilde{u}+(1-t)\mu'\nabla \tilde{u})dt\nabla w,\nabla \phi\right\ran dx\\
		=\lambda_m(\Omega)(m-1)\int_\Omega\int_{0}^{1}(t\tilde u+(1-t)\mu'\tilde u)^{m-2}w\phi dtdx,
	\end{split}
\end{equation}
i.e., 
\begin{equation}
	\nonumber
		\int_{\Omega}  \lan D^2\mathcal{H}(\nabla \tilde{u})\nabla w,\nabla \phi\ran dx
		=\lambda_m(\Omega)(m-1)\int_\Omega\tilde u^{m-2}w\phi dtdx \quad \forall  \phi\in C^\infty_0(\Omega).
\end{equation}
In order to pass the limit in $(\ref{ee7})$, we observe that for every $n,k\in \mathbb{N}$
\begin{equation}
	\nonumber
	\begin{split}
		\int_\Omega\lan A_k\nabla w_k,\nabla w_k\ran dx&\ge\int_{\{|\nabla w|\le n\}}\lan A_k\nabla w_k,\nabla w_k\ran dx\\
		&\ge\int_{\{|\nabla w|\le n\}}\lan A_k\nabla w,\nabla w\ran dx+2\int_{\{|\nabla w|\le n\}}\lan A_k\nabla w,\nabla w_k-\nabla w\ran dx.
	\end{split}
\end{equation}
Since $A_k\rightarrow D^2\mathcal{H}(t\nabla \tilde{u}+(1-t)\mu'\nabla \tilde{u})$ uniformly and $w_k\rightarrow w$ weakly in $W^{1,\theta}(\Omega)$. We get 
\begin{equation}
	\nonumber
	\lim_{k\rightarrow +\infty}\int_{\{|\nabla w|\le n\}}\lan A_k\nabla w,\nabla w_k-\nabla w\ran dx=0.
\end{equation}
This implies that for every $k\in \mathbb{N}$ we have
\begin{equation}
	\nonumber
	\begin{split}
		\liminf_{k\rightarrow \infty}\int_\Omega\lan A_k\nabla w_k,\nabla w_k\ran dx&\ge\liminf_{k\rightarrow \infty}\int_{\{|\nabla w|\le n\}}\lan A_k\nabla w,\nabla w\ran dx\\
		&=\int_{\{|\nabla w|\le n\}}\left\lan \int_{0}^{1}D^2\mathcal{H}(t\nabla \tilde{u}+(1-t)\mu'\nabla \tilde{u})dt\nabla w,\nabla w\right\ran dx\\
		&=\int_{0}^{1}(t+(1-t)\mu')^{m-2}dt\int_{\{|\nabla w|\le n\}}\lan D^2\mathcal{H}(\nabla \tilde{u})\nabla w,\nabla w\ran dx.
	\end{split}
\end{equation}
Taking  $n\rightarrow \infty$, one gets
\begin{equation}
	\label{ee10}
	\int_\Omega \lan D^2\mathcal{H}(\nabla \tilde{u})\nabla w,\nabla w\ran dx\le (m-1)\lambda_m(\Omega)\int_\Omega\tilde{u}^{m-2}w^2dx.
\end{equation}
Observe that in the right-hand side we used the strong convergence in $L^2(\Omega)$ of $\{w_k\}_{k\in \mathbb{N}}$. By recalling that 
\begin{equation}
	\nonumber
	\lan D^2\mathcal{H}(z)\xi,\xi\ran \ge|z|^{m-2}|\xi|^2, \quad \text{for every} \quad z,\xi \in \mathbb{R}^N,
\end{equation}
the estimate $(\ref{ee10})$ shows that $w$ also belongs to the weighted Sobolev space $X^{1,2}(\Omega;|\nabla u|^{m-2})$. Note also that the strong convergence of $w_k$ in $L^2(\Omega)$ together with $(\ref{ee8})$ implies that $||w||_{L^2(\Omega)}=1$, so that $w$ is non-trivial.

Finally, from the properties above we have
\begin{equation}
	\nonumber
	w\in X^{1,2}(\Omega;|\nabla \tilde u|^{m-2})\cap W^{1,1}_0(\Omega)= X^{1,2}_0(\Omega;|\nabla \tilde u|^{m-2}),
\end{equation}
\emph{Par 3.} Conclusion. From the fact that $w\in X^{1,2}_0(\Omega;|\nabla \tilde u|^{m-2})$ is nontrivial together with Proposition $\ref{pro3.2}$, Proposition $\ref{prop3.5}$ and $(\ref{ee10})$, it follows that $w$ must be proprotional either to $\tilde{u}$ or to $-\tilde{u}$. In particular, $w$ doesn't change sign.

On the other hand, by Lemma $\ref{lem4}$, we konw that $\tilde{u}_k-\tilde{v}_k$ must change sign. Accordingly, if $w^\pm_k$ stand for the positive and negative part of $w_k$ respectively, we have that each
\begin{equation}
	\nonumber
	\Omega^\pm_k:=\{x\in\Omega:w^\pm_k(x)>0\},
\end{equation}
has positive measure. Testing equation $(\ref{ee6})$ with $w^\pm_k$, we obtain by using $(\ref{ee8})$
\begin{equation}
	\label{ee11}
	\begin{split}
			\int_\Omega \lan A_k\nabla w^\pm_k,\nabla w^\pm_k\ran  dx&= (m-1)\lambda \int_\Omega\int_{0}^{1}(t\tilde{u}_k+(1-t)\tilde{v}_k)^{m-2}|w_k^\pm|^2 dtdx\\
		&+M_k^{p_k-m}(p_k-1)\int_\Omega\int_{0}^{1}(t\tilde{u}_k+(1-t)\tilde{v}_k)^{p_k-2}|w_k^\pm|^2 dtdx\\
		&\le C\int_{\Omega}|w_k^\pm|^2dx.
	\end{split}
\end{equation}
By H\"older's inequality, Theorem $\ref{thm6}$, equations $(\ref{ee9})$ and $(\ref{ee11})$ we have for an exponent $2<\sigma<\sigma_0$ 
\begin{equation}
	\nonumber
	\begin{split}
		\int_{\Omega}|w_k^\pm|^2dx&\le (\int_{\Omega}|w^\pm_k|^\sigma dx)^{\frac{2}{\sigma}} |\Omega^\pm_k|^{\frac{\sigma-2}{\sigma}}\\
		&\le \frac{1}{\mathcal{T}}|\Omega^\pm_k|^{\frac{\sigma-2}{\sigma}}\int_{\Omega}|\nabla \tilde{u}_k|^{m-2}|\nabla w^\pm_k|^2dx\\
		&\le \frac{4^{m-1}}{\mathcal{T}}|\Omega^\pm_k|^{\frac{\sigma-2}{\sigma}}\int_{\Omega}\lan A_k \nabla w^\pm_k,\nabla w^\pm_k \ran dx\\
		&\le C \frac{4^{m-1}}{\mathcal{T}}|\Omega^\pm_k|^{\frac{\sigma-2}{\sigma}}\int_{\Omega}|w_k^\pm|^2dx.
	\end{split}
\end{equation}
This implies 
\begin{equation}
	\nonumber
	|\Omega_k^\pm|\ge \frac{1}{\tilde{C}}, \quad \text{for every} \quad k\in \mathbb{N},
\end{equation}
for some constant $\tilde{C}$, not depending on $k$. This contradicts the fact that $w_k$ strongly converges in $L^2(\Omega)$ to the constant sign function $w$. The proof is over. \hfill $ \Box$
\section{Proof of the Theorems 1.1.}
\begin{lemma}
	\label{lem5.1} For $2<p<\frac{2N}{N-2}$, let $u\in W^{1,m}_0(\Omega)$ be a positive solution of equation $(\ref{eq2})$. Let $M=||u||_{L^\infty(\Omega)}$. Then there exists a constant $ \delta_1=\delta_1(N,\Omega,\lambda)>0$ such that $M\le C$ for all $m\in (2,2+\delta_1]$, where  $C=C(N,\Omega,\lambda,\delta_1)>0$.
\end{lemma}
\noindent \emph{Proof.} We argue by contradiction.  Choose $ \{m_k\}_{k\in \mathbb{N}}$ such that $\lim_{k\rightarrow +\infty}m_k=2$, and let $\{u_k\}_{k\in\mathbb{N}}\in W^{1,m}_0(\Omega)$ be the solution of $(\ref{eq2})$ with $m=m_k$ and such that $M_k:=||u_k||_{L^\infty(\Omega)}\rightarrow+\infty$. Let $x_k$ be the point where the maximum of $u_k$ is achieved. Define 
\[\tilde{u}_k(y):=\frac{u_k(\mu_ky+x_k)}{M_k},
\]
then $\tilde{u}_k$ is a function satisfying $0\leq\tilde{u}_k\leq1$, $\tilde{u}_k(0)=1$ and 
\begin{equation}
	\label{eqf13}
	\left\{
	\begin{split}
		-\Delta_{m_k} \tilde{u}_k&=\lambda \mu^{m_k}_k\tilde{u}^{{m_k}-1}_k+M_{k}^{p-m_k}\mu^{m_k}_k\tilde{u}^{p-1}_k:=f_k  \quad \text{in} \quad \Omega_k,\\
		\tilde{u}_k&=0 \quad \text{on} \quad \partial \Omega_k,
	\end{split}
	\right.
\end{equation}
where $\Omega_k=\{y\in \mathbb{R}^N:\mu_ky+x_k\in \Omega\}. $
Let $M_{k}^{p-{m_k}}\mu^{m_k}_k=1$, then $\mu_k\rightarrow0 (\text{as} \quad k\rightarrow+\infty)$. $\tilde{u}_k$ satisfies
\begin{equation}
	\label{eqf14}
	\left\{
	\begin{split}
		-\Delta_{m_k} \tilde{u}_k&=\lambda \mu^{m_k}_k\tilde{u}^{{m_k}-1}_k+\tilde{u}^{p-1}_k  \quad \text{in} \quad \Omega_k,\\
		\tilde{u}_k&=0 \quad \text{on} \quad \partial \Omega_k.
	\end{split}
	\right.
\end{equation}
Up to a subsequence, two situations may occur:
\[\text{either} \quad \text{dist}(x_k,\partial \Omega)\mu^{-1}_k\rightarrow +\infty \quad \text{or} \quad \text{dist}(x_k,\partial \Omega)\mu^{-1}_k\rightarrow d \geq 0.
\]
Case 1: $\lim_{k\rightarrow +\infty}\text{dist}(x_k,\partial \Omega)\mu^{-1}_k\rightarrow +\infty$. Then $\Omega_k\rightarrow \mathbb{R}^N$ as $k\rightarrow +\infty$.
For any $r>1$ with $B_r(0)\subset\Omega_k$, $f_k$ satisfies $(\ref{eq9})$. By Theorem $\ref{thm3}$, there exist an $\alpha\in (0,1)$ and a constant $C>0$ which doesn't depend on $k$  such that 
\[||\tilde{u}_k||_{C^{1,\alpha}(B_r(0))}\le C.
\]
By choosing a subsequence, $\tilde{u}_k\rightarrow \tilde{u}$ in $C^{1,\beta}_{loc}(\mathbb{R}^N)$ for some $\beta\in (0,\alpha)$.
Moreover, $v(0)=1$, $v\ge 0$ in $\mathbb{R}^N$, $f_k \rightarrow \tilde{u}^{p-1}$ in $C^{1,\beta}_{loc}(\mathbb{R}^N)$, and 
\[\int_{\mathbb{R}^N}\nabla \tilde{u} \nabla \psi=\int_{\mathbb{R}^N}\tilde{u}^{p-1}\psi \quad \text{for all} \quad \psi \in C^{\infty}_c(\mathbb{R}^N),
\]
namely, $\tilde{u}$ is a distributional solution of
\begin{equation}
	\label{eqf16}
	-\Delta \tilde{u}=\tilde{u}^{p-1} \quad \text{in} \quad \mathbb{R}^N.
\end{equation} 
Standard elliptic regularity yields that $\tilde{u}$ is classical solution of $(\ref{eqf16})$ which is positive. But this contradicts the classical Liouville Theorem in the whole space (see \cite[Theorem 1.2]{gidas1981priori}).

Case 2: $\text{dist}(x^k,\partial \Omega)\mu_k^{-1}\rightarrow d\ge 0$. Then $x^k\rightarrow x^0\in \partial \Omega$. Then a similar argument as in Lemma $\ref{lem2}$, we choose a subsequence, $v_k\rightarrow v \quad \text{in} \quad C^{1,\beta}_c(\mathbb{R}^N_+)$, with $||v||_{L^{\infty}(\Omega)}=1$.
$v$ satisfies 
\begin{equation}
		\label{eqf17}
		\left \{
	\begin{split}
		-\Delta{v}&={v}^{p-1} \quad  \text{in} \quad \mathbb{R}^N_+,\\
		v&=0 \quad  \text{on} \quad \partial \mathbb{R}^N_+.\\
	\end{split}
\right.
\end{equation}

By Liouville Theorem in \cite[Theorem 1.3]{gidas1981priori}, $v\equiv0$. But this contradicts $||v||_{L^{\infty}}=1$.
The lemma is proved. \hfill $ \Box$

\noindent \emph{Proof of Theorem 1.1.}  We argue by contradiction. We suppose that for every $m>2$, the equation $(\ref{eq2})$ always admits (at least) two distinct positive solutions. We then take a sequence $\{m_k\}_{k\in \mathbb{N}}\subset [2,2+\delta_1]$ such that 
\[\lim_{k\rightarrow +\infty}m_k=2.
\]
Correspondingly, for every $k\in\mathbb{N}$ there exist two distinct positive solutions $u_k$ and $v_k$ of $(\ref{eq2})$. By Theorem $\ref{thm5}$ and assumptions of Theorem $\ref{thm1}$, $u_k, v_k \rightarrow u \in C^{1,\beta}(\overline{\Omega})$, up to a subsequence. $u$ solves
\begin{equation}
	\label{eee3}
	\left\{
	\begin{split}
		-\Delta u&=\lambda u+u^{p-1} \quad \text{in} \quad  \Omega,\\
		u&=0 \quad \text{on} \quad \p \Omega.
	\end{split}
	\right.
\end{equation} 
Claim: $u\not\equiv0$. By contradiction, assume $\lim_{k\rightarrow +\infty}M_k=0$ where $M_k:=||u_k||_{L^{\infty}(\Omega)}$. Define $\tilde{u}_k=\frac{u_k}{M_k}$. Then a similar argument as in Lemma $\ref{lem3}$. $\tilde{u}_k \rightarrow \tilde{u}$ in $C^{1,\beta}(\Omega)$ where $\tilde{u}$ is the first eigenfunction of $-\Delta$. But this contradicts $\lambda < \lambda_2(\Omega)$.

The equations solved by $u_k$ and $v_k$ can be written as 
\begin{equation}
	\label{eee1}
	\int_{\Omega}|\nabla u_k|^{m_k-2}\nabla u_k\cdot \nabla \phi=\lambda \int_\Omega u_k^{m_k-1}\phi dx+\int_\Omega u^{p-1}_k\phi dx,
\end{equation} 
and
\begin{equation}
	\label{eee2}
	\int_{\Omega}|\nabla v_k|^{m_k-2}\nabla v_k\cdot \nabla \phi=\lambda \int_\Omega v_k^{m_k-1}\phi dx+\int_\Omega v^{p-1}_k\phi dx,
\end{equation} 
for any $\phi\in W^{1,m_k}_0(\Omega)$.
Let $w_k=\frac{u_k-v_k}{||u_k- v_k||_{L^2(\Omega)}}$. By using a similar way in Section 4, we get
\begin{equation}
	\label{eee5}
	\begin{split}
		\int_{\Omega}\lan A_k(x)\nabla w_k,\nabla \phi \ran dx&=\lambda (m_k-1) \int_{\Omega}\int_{0}^{1}(tu_k+(1-t)v_k)^{m_k-2}w_k\phi dtdx\\
		&+(p-1)\int_{\Omega}\int_{0}^{1}(tu_k+(1-t)v_k)^{p-2}w_k\phi dtdx,
	\end{split}
\end{equation}
where
\[A_k(x)=\int_{0}^{1}D^2\mathcal{H}_k(t\nabla u_k(x)+(1-t)\nabla v_k(x))dt,
\] 
 \[\mathcal{H}_k(z)=|z|^{m_k-2}\text{Id}+(m_k-2)|z|^{m_k-4}z\otimes z \quad \text{for} \quad  z\in \mathbb{R}^N.
 \]
Then a similar argument as in the proof of Theorem $\ref{thm2}$, $\{w_k\}_{k\in \mathbb{N}}$ converges strongly in $L^2(\Omega)$ and weakly in $W^{1,\theta}(\Omega)$ for some $\theta>1$.
In Remark $A.1$, $|\{x\in \Omega:|\nabla u(x)|=0\}|=0$. Then $\forall \epsilon _1>0$, there exists $\epsilon_2>0$ such that $|E|<\epsilon_1$ where $E=\{x\in\Omega||\nabla u|\le \epsilon_2\}$. 
\begin{equation}
	\nonumber
	\begin{split}
		|\int_{\Omega}\lan A_k(x)\nabla w_k ,\nabla \phi \ran dx-\int_\Omega\lan \nabla w,\nabla \phi\ran|&\le \int_{E^c}|\lan A_k(x)\nabla w_k,\nabla \phi\ran-\lan \nabla w,\nabla \phi\ran|dx\\
		&+\bigg|\int_{E}\lan A_k(x)\nabla w_k,\nabla \phi \ran dx\bigg|+\bigg|\int_{E}\lan \nabla w,\nabla \phi \ran dx\bigg|
	\end{split}
\end{equation}
The last two termconverges to 0 by the continuity of integral.
Thus $w$ satisfies 
\begin{equation}
	\nonumber
	\int_{\Omega}\lan \nabla w,\nabla \phi \ran dx=\lambda \int_{\Omega} w\phi +(p-1)\int_{\Omega}u^{p-2}w\phi dx.
\end{equation}
This contradicts the nondegeneracy of the equation $(\ref{eq04})$. \hfill $\Box $

\section*{ Appendix A}
\noindent \textbf{A.1. The linearized equation. } We first observe that $u\in C^{1,\beta}(\Omega)$ and the critical set 
\begin{equation}
	Z:=\{x\in \Omega :|\nabla u|=0 \} \nonumber
\end{equation}
is a compact set contained in $\Omega$, thanks to Theorem $\ref{thm5}$. Thus $\Omega \backslash Z$ is an open set. On this set, by classical Elliptic Regularity we can infer that $u \in C^2(\Omega \backslash Z)$. 

\noindent \textbf{Remark A.1} (Hessian terms)\textbf{.} We seize the opportunity to mention that, since $u\in L^{\infty}(\Omega)$. The right-hand side  of $(\ref{eq27})$ is bounded. Then we have 
\begin{equation}
	\tag{A.0}
	|\nabla u|^{m-2}\nabla u\in W^{1,2}_{\mathrm{loc}}(\Omega),
\end{equation} 
thanks to \cite[Theorem 2.1]{cianchi2018second}. In addition, as noted in \cite[Remark 2.3]{damascelli2004}, the weak gradient of $|\nabla u|^{m-2}\nabla u$ coincides with the classical gradient in $\Omega \backslash Z$, while
\begin{equation}
	\nabla (|\nabla u|^{m-2}\nabla u)=0, \quad \text{a.e. on} \quad Z,\nonumber
\end{equation}
since by definition $Z$ coincides with the zero level set of $|\nabla u|^{m-2}\nabla u$, thus it is sufficient to use a standard property of Sobolev functions. By futher observing that $|Z|=0$ (by virtue of \cite[Theorem 2.3]{damascelli2004}). We can actually say that the weak gradient of $|\nabla u|^{m-2}\nabla u$ coincides with the classical gradient almost everywhere in $\Omega$.

We consider the following equation
\begin{equation}
	\tag{A.1}
	\left\{
	\begin{split}
		-\Delta_m u&=f(u) \quad \text{in} \quad  \Omega,\\
		u&=0 \quad \text{on} \quad \partial \Omega,
	\end{split}
	\right.
\end{equation}
where  $f\in C^{1,1}(\mathbb{R}^+)$. 
We assume that $u$ satisfies the following properties.

\noindent (u1)  There exist $\chi\in (0,1)$ and $L>0$ such that 
\[||u||_{C^{1,\chi}(\Omega)}\le L.
\]
\noindent  (u2) There exist $\tau>0$, $\mu_0>0$ and $\mu_1>0$, such that 
\[|\nabla u|\ge \mu_0, \quad \text{in} \quad \Omega_{\tau}, 
\]
and 
\[u\ge \mu_1, \quad \text{in} \quad \overline{\Omega \backslash \Omega_{\tau}},
\]
where $\Omega_{\tau}= \{x\in \Omega: \text{dist}(x,\p\Omega)\le \tau\}$.

We then take $\psi\in C^{\infty}_0(\Omega \backslash Z)$ and test the weak formulation of $(A.1)$ against a partial derivative $\psi_{x_i}$. We obtain 
\begin{equation}
	\tag{A.2}
	\begin{split}
		\int_{\Omega}\lan |\nabla u|^{m-2}\nabla u_{x_i},\nabla \psi\ran dx&+(m-2)\int_{\Omega}|\nabla u|^{m-4}\lan \nabla u,\nabla u_{x_i} \ran \lan \nabla u,\nabla \psi \ran dx\\
		&=\int_{\Omega} f'(u)u_{x_i}\psi dx  \quad \text{for every}
		\quad \psi \in C^{\infty}_0(\Omega \backslash Z). \nonumber 
	\end{split}
\end{equation}

\noindent \textbf{Proposition A.2.} Let $m>2$ and let $\Omega\subset \mathbb{R}^N$ be an open bounded connected set, with boundary of class $C^{1,\alpha}$, for some $0<\alpha<1$. Let $u\in W^{1,m}_0(\Omega)$ be a positive solution of $(A.1)$. Assume $u$ satisfies (u1) and (u2). Let $\beta \in [0,1)$ and 
\begin{equation}
	\left\{
	\begin{split}
		\gamma <N-2, \quad \text{if} \quad  N\ge 3,\\
		\gamma \le 0, \quad \text{if} \quad N=2. \nonumber
	\end{split}
	\right.
\end{equation}
Then for every $i\in\{1,\cdots, N\}$, if we set $Z_i=\{y\in \Omega:u_{x_i}(y)=0\}$, we have the estimate 
\begin{equation}
	\tag{A.3}
	\sup_{x\in \Omega}\int_{E \backslash Z_i} \frac{|\nabla u(y)|^{m-2}|u_{x_i}(y)|^{-\beta}|\nabla u_{x_i}(y)|^2}{|x-y|^{\gamma}}dy \le C_1
\end{equation}
for some $C_1=C_1(N, L, ||f||_{C^{1,1}(\mathbb{R})}, \alpha, \beta, \gamma, \Omega, \text{dist}(E,\p\Omega))(1+(m-1)L^{m-\beta})>0$.

\noindent \emph{Proof.} Without loss of generality, we prove the result for $\gamma \geq 0$. The heuristic idea is to test the linearized equation $(A.2)$ with $u_{x_i}\vline u_{x_i}\vline^{-\beta}|x-y|^{-\gamma} \phi ^2$, where $\phi$ is a smooth cut-off function. In order to do this rigorously, we fix $x\in \Omega$ and for $0<\epsilon <1$ use the test function 
\begin{equation}
	\psi(y)=G_{\epsilon}(u_{x_i}(y))|u_{x_i}(y)|^{-\beta}(|x-y|+\epsilon)^{-\gamma} \phi ^2(y),\nonumber
\end{equation} 
where $\phi\in C^{\infty}_0(\Omega)$ is such that
\begin{equation}
	\phi\equiv 1\quad \text{on} \quad E, \quad 0\le \phi \le1, \quad ||\nabla\phi||_{L^{\infty}(E)}\le \frac{C}{\text{dist}(E,\partial \Omega)},\nonumber
\end{equation}
and the odd Lipschitz function $G_{\epsilon}$ is given by 
\begin{equation}
	G_{\epsilon}(t)=\max \{t-\epsilon,0\}, \quad \text{for}\quad  t\ge 0, \quad \text{and} \quad G_{\epsilon}(t)=-G_{\epsilon}(-t), \quad \text{for} \quad t\le 0 \nonumber 
\end{equation}
The function $\psi$ is a product of a Lipschitz function of $u_{x_i}$, which vanishes in a neighborhood of the critical set $Z$, and a smooth function with compact support in $\Omega$. Therefore $\psi$ is an admissible test function for $(A.2)$. This gives
\begin{equation}
	\begin{split}
		&\int_{\Omega}\frac{|\nabla u|^{m-2}|\nabla u_{x_i}|^2|u_{x_i}|^{-\beta}}{(|x-y|+\epsilon)^{\gamma}}\left[G'_{\epsilon}(u_{x_i})-\beta \frac{G_{\epsilon}(u_{x_i})}{u_{x_i}}\right]\phi ^2dy\\
		&+(m-2)\int_{\Omega}\frac{|\nabla u|^{m-4}(\lan \nabla u,\nabla u_{x_i}\ran)^2|u_{x_i}|^{-\beta}}{(|x-y|+\epsilon)^{\gamma}}\left[G'_{\epsilon}(u_{x_i})-\beta \frac{G_{\epsilon}(u_{x_i})}{u_{x_i}}\right]\phi ^2dy\\
		&+2\int_{\Omega}\frac{|\nabla u|^{m-2}\lan \nabla u_{x_i},\nabla \phi \ran|u_{x_i}|^{-\beta}}{(|x-y|+\epsilon)^{\gamma}}\phi G_{\epsilon}(u_{x_i})dy\\
		&+2(m-2)\int_{\Omega}\frac{|\nabla u|^{m-4}\lan \nabla u,\nabla u_{x_i}\ran \lan \nabla u,\nabla \phi \ran|u_{x_i}|^{-\beta}}{(|x-y|+\epsilon)^{\gamma}}\phi G_{\epsilon}(u_{x_i})dy\\
		&+\int_{\Omega}|\nabla u|^{m-2}\lan \nabla u_{x_i},\nabla ((|x-y|+\epsilon)^{-\gamma})\ran G_{\epsilon}(u_{x_i})|u_{x_i}|^{-\beta}\phi ^2dy\\
		&+(m-2)\int_{\Omega}|\nabla u|^{m-4}\lan \nabla u,\nabla u_{x_i}\ran \left\lan \nabla u,\nabla ((|x-y|+\epsilon)^{-\gamma})\right\ran G_{\epsilon}(u_{x_i})|u_{x_i}|^{-\beta}\phi ^2dy\\
		&= \int_{\Omega}\frac{f'(u)G_{\epsilon}(u_{x_i})|u_{x_i}|^{-\beta}}{(|x-y|+\epsilon)^{\gamma}}\phi ^2dy\\
	\nonumber
    \end{split}
\end{equation}
Note that 
\begin{equation}
	G'_{\epsilon}(t)-\beta\frac{G_{\epsilon}(t)}{t}\ge 0, \quad \text{for every} \quad t\in \mathbb{R},\nonumber
\end{equation}
therefore, by dropping the second term on the left-hand side and using the Cauchy-Schwarz inequality, we get 
\begin{equation}
	\begin{split}
		&\int_{\Omega}\frac{|\nabla u|^{m-2}|\nabla u_{x_i}|^2|u_{x_i}|^{-\beta}}{(|x-y|+\epsilon)^{\gamma}}\left[G'_{\epsilon}(u_{x_i})-\beta \frac{G_{\epsilon}(u_{x_i})}{u_{x_i}}\right]\phi ^2dy\\
		&\le 2(m-1)\int_{\Omega}\frac{|\nabla u|^{m-2} |\nabla u_{x_i}| |u_{x_i}|^{-\beta}}{(|x-y|+\epsilon)^{\gamma}} G_{\epsilon}(u_{x_i})\phi|\nabla \phi|dy\\
		&+\gamma (m-1)\int_{\Omega}\frac{|\nabla u|^{m-2} |\nabla u_{x_i}| |u_{x_i}|^{-\beta}}{(|x-y|+\epsilon)^{\gamma+1}} G_{\epsilon}(u_{x_i})\phi^2dy\\
		&+ \int_{\Omega}\frac{ f'(u) G_{\epsilon}(u_{x_i})  |u_{x_i}|^{-\beta}}{(|x-y|+\epsilon)^{\gamma}} \phi^2dy\\
		&:=I_1+I_2+I_3.
	\end{split}
	\tag{A.4}
\end{equation}
By using Theorem $\ref{thm5}$, the properties of the cut-off function $\phi$, the last term can be estimated as 
\begin{equation}
	I_3\le C\int_{\Omega}\frac{1}{|x-y|^{\gamma}}dy, \nonumber
\end{equation}
for a constant $C$ depending on $N$, $||f||_{C^{1,1}(\mathbb{R})}$, $L$, $\alpha$, $\beta$ and $\Omega$, only. In turn, the last integral is easily estimated as follows 
\begin{equation}\tag{A.5}
	\begin{split}
		\int_{\Omega}\frac{1}{|x-y|^{\gamma}}dy&\le \int_{\{y\in \mathbb{R}^N:|y-x|\le \text{diam}(\Omega)\}}\frac{1}{|y-x|^{\gamma}}dy\\
		&=N\omega_{N}\int_{0}^{\text{diam}(\Omega)}\rho ^{N-1-\gamma}d\rho =\frac{N\omega_N}{N-\gamma}\left(\text{diam}(\Omega)\right)^{N-\gamma}.
	\end{split}
\end{equation}
In the last integral we used that $\gamma<N-2$.

As for the terms $I_1$ and $I_2$, we first observe that by using 
\begin{equation}
	\tag{A.6}
	(|x-y|+\epsilon)^{\gamma}\ge \frac{(|x-y|+\epsilon)^{\gamma+1}}{\text{diam}(\Omega)+1}, \quad \text{for every} \quad x,y\in \Omega,
\end{equation}
we have 
\begin{equation}
	\nonumber
	I_1+I_2\le (2+|\gamma|)(m-1)\int_{\Omega}\frac{|\nabla u|^{m-2} |\nabla u_{x_i}| |u_{x_i}|^{-\beta}}{(|x-y|+\epsilon)^{\gamma+1}}G_{\epsilon}(u_{x_i})\phi [\phi +|\nabla \phi|]dy.
\end{equation}
 Thus, from $(A.4)$, we have obtained 
\begin{equation}
	\tag{A.7}
	\begin{split}
		\int_{\Omega}&\frac{|\nabla u|^{m-2} |\nabla u_{x_i}|^2 |u_{x_i}|^{-\beta}}{(|x-y|+\epsilon)^{\gamma}} \left[G'_{\epsilon}(u_{x_i})-\beta\frac{G_{\epsilon}(u_{x_i})}{u_{x_i}}\right]\phi^2dy\\
		&\le C+(2+|\gamma|)(m-1)\int_{\Omega}\frac{|\nabla u|^{m-2} |\nabla u_{x_i}| |u_{x_i}|^{-\beta}}{(|x-y|+\epsilon)^{\gamma+1}}G_{\epsilon}(u_{x_i})\phi [\phi +|\nabla \phi|]dy.
	\end{split}
\end{equation}
In order to estimate the last integral, we used Young's inequality. For every $\delta>0$, we have 
\begin{equation}
	\tag{A.8}
	\begin{split}
		(2+|\gamma|)(m-1)&\int_{\Omega}\frac{|\nabla u|^{m-2} |\nabla u_{x_i}| |u_{x_i}|^{-\beta}}{(|x-y|+\epsilon)^{\gamma+1}}G_{\epsilon}(u_{x_i})\phi [\phi +|\nabla \phi|]dy\\
		&\le \frac{\delta}{2} \int_{\Omega} \frac{|\nabla u|^{m-2} |\nabla u_{x_i}|^2 |u_{x_i}|^{-\beta}}{(|x-y|+\epsilon)^{\gamma}} \frac{G_{\epsilon}(u_{x_i})}{u_{x_i}}\phi^2dy\\
		&+\frac{((2+|\gamma|)(m-1))^2}{2\delta}  \int_{\Omega} \frac{|\nabla u|^{m-2}  |u_{x_i}|^{1-\beta}}{(|x-y|+\epsilon)^{\gamma+2}}G_{\epsilon}(u_{x_i}) [\phi+|\nabla \phi|]^2dy.
	\end{split}
\end{equation}
We make the choice $\delta=(2+|\gamma|)(m-1)$. Then observe that
\begin{equation}
       \nonumber
       (1-\beta)\frac{G_{\epsilon}}{t}\le G'_{\epsilon}(t)-\beta\frac{G_{\epsilon}(t)}{t},\quad \text{for  every} \quad |t|\neq \epsilon.
\end{equation}
Then the first term in the right-hand side of $(A.8)$ can now be estimated by
\begin{equation}
	\nonumber
	\frac{(2+|\gamma|)(m-1)}{2}\int_{\Omega}\frac{|\nabla u|^{m-2} |\nabla u_{x_i}|^2 |u_{x_i}|^{-\beta}}{(|x-y|+\epsilon)^{\gamma}}\left[G'_{\epsilon}(u_{x_i})-\beta \frac{G_{\epsilon}(u_{x_i})}{u_{x_i}}\right]\phi^2dy,
\end{equation}
which can be absorbed into the right-hand side of $(A.7)$. Thus, up to now, we obtained 
\begin{equation}
	\tag{A.9}
	\begin{split}
		\int_{\Omega}&\frac{|\nabla u|^{m-2} |\nabla u_{x_i}|^2 |u_{x_i}|^{-\beta}}{(|x-y|+\epsilon)^{\gamma}}\left[G'_{\epsilon}(u_{x_i})-\beta \frac{G_{\epsilon}(u_{x_i})}{u_{x_i}}\right]\phi^2dy\\
		&\le C+(2+|\gamma|)(m-1)\int_{\Omega} \frac{|\nabla u|^{m-2}  |u_{x_i}|^{1-\beta}}{(|x-y|+\epsilon)^{\gamma+2}}G_{\epsilon}(u_{x_i}) [\phi+|\nabla \phi|]^2dy,
	\end{split}
\end{equation}
possibly for a different constant, independent of $x\in \Omega$, and $\epsilon \in (0,1)$. Using that $|G_{\epsilon}(t)|\le |t|$ and the properties of $\phi$, the last integral of $(A.9)$ can be estimated by
\begin{equation}
	\nonumber
	\int_{\Omega} \frac{|\nabla u|^{m-2}  |u_{x_i}|^{2-\beta}}{(|x-y|+\epsilon)^{\gamma+2}} [\phi+|\nabla \phi|]^2dy\le CL^{m-\beta}\int_{\Omega}\frac{1}{|x-y|^{\gamma +2}}dy,
\end{equation}
and the last integral is uniformly (in $x\in \Omega$) bounded by a constant depending only on $N$ and $\text{diam}(\Omega)$.
 From $(A.9)$ and using that $\phi\equiv 1$ on $E$, we thus obtain 
 \begin{equation}
 	\tag{A.10}
 	\int_{E}\frac{|\nabla u|^{m-2} |\nabla u_{x_i}|^2 |u_{x_i}|^{-\beta}}{(|x-y|+\epsilon)^{\gamma}}\left[G'_{\epsilon}(u_{x_i})-\beta \frac{G_{\epsilon}(u_{x_i})}{u_{x_i}}\right]\phi^2dy\le C(1+(m-1)L^{m-\beta}),
 \end{equation}
for some $C=C(N, L, ||f||_{C^{1,1}(\mathbb{R})} , \gamma, \alpha, \beta, \Omega, \text{dist}(E,\p \Omega) )>0$.  We have
\begin{equation}
	\tag{A.11}
	\lim_{\epsilon \rightarrow 0}\left[G'_{\epsilon}-\beta \frac{G_{\epsilon}(\tau)}{\tau}\right]=1-\beta, \quad \text{for} \quad \tau \neq 0.
\end{equation}
By using Fatou's Lemma together with $(A.11)$, we may take the limit as $\epsilon$ goes to $0$ in $(A.10)$ and obtain 
\begin{equation}
	\nonumber
	\int_{E\backslash Z_i}\frac{|\nabla u|^{m-2}|u_{x_i}|^{-\beta}|\nabla u_{x_i}|^2}{|x-y|^{\gamma}}dy\le C.
\end{equation}
The previous bound holds uniformly with respect to $x\in \Omega$, thus we get the desired conclusion. \hfill $ \Box$

\noindent \textbf{Corollary A.3.} Under the assumption of Proposition A.2, for every $\beta \in (-\infty,1)$, we have 
\begin{equation}
	\nonumber
	\begin{split}
		\sup_{x\in \Omega}\int_{E}\frac{|\nabla u(y)|^{m-2-\beta}|D^2u(y)|^2}{|x-y|^{\gamma}}dy\le C_1\left(1+(m-1)L^{m-\beta}\right).
	\end{split}
\end{equation}
Here we still denote by $C_1$ the same constant as in Proposition A.2.

\noindent \emph{Proof.} From (A.3), by using that $|u_{x_i}|^{-\beta}\ge |\nabla u|^{-\beta}$, we immediately get
\begin{equation}
	\nonumber
	\sup_{x \in \Omega}\int_{E\backslash Z_i}\frac{|\nabla u(y)|^{m-2-\beta}|\nabla u_{x_i}(y)|^2}{|x-y|^{\gamma}}dy\le C_1\left(1+(m-1)L^{m-\beta}\right).
\end{equation}
We then observe that $u_{x_i}\in C^1(\Omega \backslash Z)$ and thus it belongs to $W^{1,1}_{\text{loc}}(\Omega \backslash Z)$. By appealing to \cite[Theorem 6.19]{lieb2001analysis}, we have 
\begin{equation}
	\nonumber
	\nabla u_{x_i}=0, \quad a.e.\quad \text{in}\quad  Z_i\cap(\Omega \backslash Z).
\end{equation}
By using this fact and summing over $i=1,\cdots, N$,  we get the claimed inequality, by recalling that $|Z|=0$, see Remark $A.1$. 

The case $\beta <0$ can be reduced to the case $\beta =0$: it is sufficient to use that $||\nabla u||_{L^{\infty}}\le L\le +\infty$, thus we get 
\[\sup_{x \in \Omega}\int_{E}\frac{|\nabla u(y)|^{m-2-\beta}|D^2u(y)|^2}{|x-y|^{\gamma}}dy\le L^{-\beta}\sup_{x \in \Omega}\int_{E}\frac{|\nabla u(y)|^{m-2}|D^2u(y)|^2}{|x-y|^{\gamma}}dy.
\]
\hfill $ \Box$

\noindent \textbf{Proposition A.4.} Let $m>2$ and let $\Omega\subset \mathbb{R}^N$ be an open bounded connected set, with boundary of class $C^{1,\alpha}$, for some $0<\alpha<1$. Let $\lambda >0$ and let $f(t)-\lambda t^{m-1}\ge 0$ for all $t$. Let $u\in W^{1,m}_0(\Omega)$ be a positive solution of $(A.1)$.  Assume $u$ satisfies (u1). Let  
\begin{equation}
	\left\{
	\begin{split}
		\gamma <N-2, \quad \text{if} \quad  N\ge 3,\\
		\gamma \le 0, \quad \text{if} \quad N=2. \nonumber
	\end{split}
	\right.
\end{equation}
Then for every $K\Subset E\Subset \Omega$ and every $b<1$, we have 
\begin{equation}
	\nonumber
	\sup_{x\in \Omega} \int_{K} \frac{1}{|\nabla u|^{(m-1)b}|x-y|^{\gamma}}dy\le C_2\left(\frac{1}{\inf_{E} u}\right)^{m-1} \left(L^{(m-1)(1-b)}+\frac{(m-1)^3L^{m-b}}{(\inf_{E}u)^{m-1}}\right),
\end{equation}
where $C_2=C_2(N, L, \alpha, b, \gamma,\lambda, \text{dist}(K,\p E), \text{dist}(E, \p  \Omega))>0$.

\noindent \emph{Proof.} Without loss of generality, we prove the result for $\gamma \ge 0$. The heurstic idea is to test equation $(A.1)$ with $|\nabla u|^{-(m-1)b}|x-y|^{-\gamma}\phi$, where $\phi\in C^{\infty}_0(E)$ is a cut-off function such that
\begin{equation}
	\nonumber
	\phi\equiv 1 \quad  \text{in}\quad  K,\quad 0\le \phi \le 1, \quad ||\phi||_{L^{\infty}(E)}\le \frac{C}{\text{dist}(K,\p E)}. 
\end{equation}
To make this precise, we fix $x\in \Omega$ and use for every $\epsilon >0$ the test function
\begin{equation}
	\tag{A.12}
	\psi(y)=(|\nabla u|^{m-1}+\epsilon)^{-b}(|x-y|+\epsilon)^{-\gamma}\phi (y).
\end{equation}
This is a product of a Lipschitz function of $|\nabla u|^{m-1}$ and a smooth function with compact support in $\Omega$. In light of $(A.2)$, we have that this function belongs to $W^{1,2}_0(\Omega)$. If we now use that $u\in L^{\infty}(\Omega)$, we see that in the weak formulation of $(A.1)$ we can in particular admit test functions $\psi \in W^{1,2}_0(\Omega)$. Therefore the test function in $(A.12)$ is feasible. 

This gives 
\begin{equation}
	\nonumber
	\begin{split}
		\lambda \int_{\Omega}\frac{u^{m-1}}{(|\nabla u|^{m-1}+\epsilon)^{b}} &\frac{\phi}{(|x-y|+\epsilon)^{\gamma}}dy\\
		+&\int_{\Omega}\frac{f(u)-\lambda u^{m-1}}{(|\nabla u|^{m-1}+\epsilon)^{b}} \frac{\phi}{(|x-y|+\epsilon)^{\gamma}}dy\\
		=&\int_{\Omega}\frac{\lan|\nabla u|^{m-2}\nabla u, \nabla \phi \ran}{(|\nabla u|^{m-1}+\epsilon)^{b}}\frac{1}{(|x-y|+\epsilon)^{\gamma}}dy\\
		-&\beta(m-1)\int_{\Omega}\frac{|\nabla u|^{2m-5}}{(|\nabla u|^{m-1}+\epsilon)^{b +1}}\lan D^2u\nabla u,\nabla u\ran\frac{\phi}{(|x-y|+\epsilon)^{\gamma}}dy\\
		+&\int_{\Omega}\frac{|\nabla u|^{m-2}}{(|\nabla u|^{m-1}+\epsilon)^{b}}	\lan \nabla u,\nabla ||x-y|+\epsilon|^{-\gamma} \ran \phi dy.
	\end{split}
\end{equation}
Dropping the second term in the left hand side and using that 
\begin{equation}
	\nonumber
	u\geq \inf_{E} u>0,
\end{equation}
we obtain
\begin{equation}
	\tag{A.13}
	\begin{split}
		(\inf_{E} u)^{m-1} \int_{\Omega}\frac{1}{(|\nabla u|^{m-1}+\epsilon)^{b}} &\frac{\phi}{(|x-y|+\epsilon)^{\gamma}}dy\\
		\le &C\int_{\Omega}\frac{|\nabla u|^{m-1} |\nabla \phi |}{(|\nabla u|^{m-1}+\epsilon)^{b}}\frac{1}{(|x-y|+\epsilon)^{\gamma}}dy\\
		+&Cb(m-1)\int_{\Omega}\frac{|\nabla u|^{2m-5}}{(|\nabla u|^{m-1}+\epsilon)^{b +1}}\lan D^2u\nabla u,\nabla u\ran\frac{\phi}{(|x-y|+\epsilon)^{\gamma}}dy\\
		+&C\gamma \int_{\Omega}\frac{|\nabla u|^{m-1}}{(|\nabla u|^{m-1}+\epsilon)^{b}}	\frac{\phi}{(|x-y|+\epsilon)^{\gamma+1}}dy\\
		:=&J_1+J_2+J_3,
	\end{split}
\end{equation}
where $C=C(N,\Omega, \lambda)>0$. We observe that, by using $(A.6)$, we have 
\begin{equation}
	\nonumber
	J_1+J_3\le C\int_{\Omega}\frac{|\nabla u|^{m-1}}{(|\nabla u|^{m-1}+\epsilon)^{b}}	\frac{\phi+|\nabla \phi|}{(|x-y|+\epsilon)^{\gamma+1}}dy,
\end{equation}
for a constant $C$ depending on $\text{diam}(\Omega)$ and $\gamma$. We can then go by observing that $|\nabla u|^{m-1}\le |\nabla u|^{m-1}+\epsilon$. This gives 
\begin{equation}
	\nonumber
	J_1+J_3\le CL^{(m-1)(1-b)}\int_{\Omega}\frac{1}{(|x-y|+\epsilon)^{\gamma+1}}dy,
\end{equation}
for a constant $C=C(N, \alpha, b, \gamma, \Omega, \text{dist}(K,\p E),\lambda)>0$. Then we can estimate the last integral as in $(A.5)$. For $J_2$ we have
\begin{equation}
	\nonumber
	J_2\le C(m-1)\int_{\Omega}\frac{|\nabla u|^{2m-3}}{(|\nabla u|^{m-1}+\epsilon)^{b+1}}	\frac{|D^2u|}{(|x-y|+\epsilon)^{\gamma}}\phi dy.
\end{equation}
In the last integral above, we use Young's inequality as follows
\begin{equation}
	\nonumber
	\begin{split}
		C(m-1)\int_{\Omega}\frac{|\nabla u|^{2m-3}}{(|\nabla u|^{m-1}+\epsilon)^{b+1}}&\frac{|D^2u|}{(|x-y|+\epsilon)^{\gamma}}\phi dy\\
		\le&\frac{(\inf_{E}u)^{m-1}}{2}\int_{\Omega}\frac{1}{(|\nabla u|^{m-1}+\epsilon)^{b}}\frac{\phi}{(|x-y|+\epsilon)^{\gamma}} dy\\
		+&\frac{(C(m-1))^2}{2(\inf_{E}u)^{m-1}} \int_{\Omega}\frac{|\nabla u|^{4m-6}}{(|\nabla u|^{m-1}+\epsilon)^{b+2}}\frac{|D^2u|^2\phi}{(|x-y|+\epsilon)^{\gamma}} dy\\
		\le & \frac{(\inf_{E}u)^{m-1}}{2}\int_{\Omega}\frac{1}{(|\nabla u|^{m-1}+\epsilon)^{b}}\frac{\phi}{(|x-y|+\epsilon)^{\gamma}} dy\\
		+&\frac{(C(m-1))^2}{2(\inf_{E}u)^{m-1}} \int_{\Omega} |\nabla u|^{(2-b)(m-2)-b}\frac{|D^2u|^2\phi}{(|x-y|+\epsilon)^{\gamma}} dy.\\
	\end{split}
\end{equation}
The last integral is uniformly bounded, thanks to Corollary $(A.3)$. We then obtain from $(A.13)$ 
\begin{equation}
	\nonumber
	\begin{split}
		(\inf_{E}u)^{m-1}&\int_{\Omega}\frac{1}{(|\nabla u|^{m-1}+\epsilon)^{\beta}}\frac{\phi}{(|x-y|+\epsilon)^{\gamma}} dy
		\le C\left(L^{(m-1)(1-b)}+\frac{(m-1)^3L^{m-b}}{(\inf_{E}u)^{m-1}}\right)\\
		&+\frac{(\inf_{E}u)^{m-1}}{2} \int_{\Omega}\frac{1}{(|\nabla u|^{m-1}+\epsilon)^{\beta}}\frac{\phi}{(|x-y|+\epsilon)^{\gamma}} dy,
	\end{split}
\end{equation}
where the constant $C=C(N, \alpha, b, \gamma, \lambda, \Omega, \text{dist}(E,,\p E))>0$. The term on the right-hand side can now be absorbed in the left-hand side. Since $\phi =1$ on $K$, this implies the desired result upon letting $\epsilon$ go to $0$ and using Fatou's Lemma.\hfill $ \Box$

\noindent \textbf{Theorem A.5.} Let $m>2$ and let $\Omega\subset \mathbb{R}^N$ be an open bounded connected set, with boundary of class $C^{1,\alpha}$, for some $0<\alpha<1$. Let $\lambda >0$ and let $f(t)-\lambda t^{m-1}\ge 0$ for all $t$. Let $u\in W^{1,m}_0(\Omega)$ be a positive solution of $(A.1)$.  Assume $u$ satisfies (u1) and (u2). Let  
\begin{equation}
	\left\{
	\begin{split}
		\gamma <N-2, \quad \text{if} \quad  N\ge 3,\\
		\gamma \le 0, \quad \text{if} \quad N=2. \nonumber
	\end{split}
	\right.
\end{equation}
Then every $r<m-1$, there exist $\mathcal{S}=\mathcal{S}(\alpha, N, \Omega, L, \mu_0, \mu_1, r, \gamma, \lambda, m)>0$ such that 
\begin{equation}
	\nonumber
	\sup_{x\in \Omega}\int_{\Omega}\frac{1}{|\nabla u_p(y)|^{r}|y-x|^{\gamma}}dy\le \mathcal{S}.
\end{equation}
In particular, we also have
\begin{equation}
	\nonumber
	\tag{A.14}
	\sup_{\Omega}\int_{\Omega}\frac{1}{|\nabla u_p(y)|^{r}}dy\le \tilde{\mathcal{S}},
\end{equation}
for some $\tilde{\mathcal{S}}=\tilde{\mathcal{S}}(\alpha, N, \Omega, L, \mu_0, \mu_1, r, \gamma, \lambda, m)>0$.

\noindent \emph{Proof.} We only treat the case $N\ge 3$ in detail. This is \cite[Theorem 2.3]{damascelli2004}: as claimed, we just want to pay attention to the dependence of the constant $\mathcal{S}$ on the data. We have 
\begin{equation}
	\tag{A.15}
	\int_{\Omega_{\tau}}\frac{1}{|\nabla u_p(y)|^r|y-x|^{\gamma}}dy\le \frac{1}{\mu_0^r}\int_{\Omega_{\tau}}\frac{1}{|y-x|^{\gamma}}dy\le \frac{N\omega_{N}}{(N-\gamma)\mu_0^r}(\text{diam}(\Omega))^{N-\gamma}.
\end{equation}
Thus we have a uniform estimate, at least when integrating in fixed neighborhood of the boundary.
In order to prove a uniform estimate on 
\begin{equation}
	\nonumber
	\int_{\Omega\backslash \Omega_{\tau}}\frac{1}{|\nabla u_p(y)|^r|y-x|^{\gamma}}dy, \quad \text{for  every} \quad x\in \Omega,
\end{equation}
we apply Proposition $A.4$ with 
\begin{equation}
	\nonumber
	E=\Omega \backslash \Omega_{\tau/2},\quad  K=\Omega\backslash \Omega_{\tau}, \quad b=\frac{r}{m-1},
\end{equation}
and with the constant $\mu_1$ provided by Theorem $\ref{thm5}$. This yields 
\begin{equation}
	\nonumber
	\int_{\Omega\backslash \Omega_{\tau}}\frac{1}{|\nabla u_p(y)|^r|y-x|^{\gamma}}dy\le C_2\left(\frac{1}{\inf_{E} u}\right)^{m-1} \left(L^{(m-1)(1-b)}+\frac{(m-1)^3L^{m-b}}{(\inf_{E}u)^{m-1}}\right),
\end{equation}
as desired.

 Finally, the estimate $(A.14)$ is an easy consequence of the previous one, it is sufficient to take $\gamma=0$.\hfill $\Box$

\bibliographystyle{plain}
\bibliography{citeliouville}

{\em Addresse and E-mail:}
\medskip

{\em Wei Ke}

{\em School of Mathematical Sciences}

{\em Fudan University}

{\em wke21@m.fudan.edu.cn}

\end{document}